\documentclass[12pt,twoside]{article}
\usepackage[hmargin=0.8in,vmargin=0.9in]{geometry}
\usepackage{fancyhdr}
\usepackage{graphicx}
\usepackage{subfigure}
\usepackage{amssymb}
\usepackage{amsmath}
\usepackage{amsfonts}
\usepackage{theorem}
\usepackage{mathrsfs}
\usepackage{mathtools}
\usepackage{bm}
\usepackage{color}
\usepackage{setspace}
\usepackage{exscale}
\usepackage{relsize}
\usepackage{epstopdf}
\usepackage{float}
\usepackage{cite}
\usepackage{hyperref}

\usepackage{latexsym}
\usepackage{empheq}
\usepackage{bbm}
\usepackage{xcolor}
\usepackage{tikz}
\usepackage{cases}
\usetikzlibrary{calc,intersections}
\usepackage{nicefrac}
\usepackage{url}
\usepackage{booktabs} 
\usepackage{array} 
\usepackage{paralist} 
\usepackage{verbatim} 
\usepackage{subfigure} 
\usepackage[ruled,vlined,algosection]{algorithm2e}

\theoremstyle{plain}			
\newtheorem{thm}{Theorem}[section]

\newtheorem{rmk}[thm]{Remark}

\newcommand*\xbar[1]{%
  \hbox{%
    \vbox{%
      \hrule height 0.5pt 
      \kern0.4ex
      \hbox{%
        \kern-0.05em
        \ensuremath{#1}%
        \kern-0.00em
      }%
    }%
  }%
}

\newcommand\eref[1]{(\ref{#1})}

\allowdisplaybreaks[1]

\numberwithin{equation}{section}
\numberwithin{figure}{section}
\numberwithin{table}{section}

\newcommand{\dt}{\Delta t}
\newcommand{\dx}{\Delta x}
\newcommand{\dy}{\Delta y}

\newcommand{\Bx}{\bm{x}}
\newcommand{\Bu}{\bm{u}}

\title{A Hybrid Finite-Difference-Particle Method for Chemotaxis Models}
\author{Alina Chertock\thanks{Department of Mathematics, North Carolina State University, Raleigh, NC 27695, USA;
{\tt chertock@math.ncsu.edu}}, Shumo Cui\thanks{Shenzhen International Center for Mathematics, Southern University of Science and Technology
(SUSTech), Shenzhen, 518055, PR China; {\tt cuism@sustech.edu.cn}}, Alexander Kurganov\thanks{Department of Mathematics, Shenzhen
International Center for Mathematics, and Guangdong Provincial Key Laboratory of Computational Science and Material Design, Southern
University of Science and Technology (SUSTech), Shenzhen 518055, PR China; {\tt alexander@sustech.edu.cn}}, and Chenxi Wang\thanks{Shenzhen
International Center for Mathematics, Southern University of Science and Technology (SUSTech), Shenzhen, 518055, PR China;
\tt{wangcx2017@mail.sustech.edu.cn}}}

\begin{document}

\date{}
\maketitle

\begin{abstract}
Chemotaxis systems play a crucial role in modeling the dynamics of bacterial and cellular behaviors, including propagation, aggregation, and
pattern formation, all under the influence of chemical signals. One notable characteristic of these systems is their ability to simulate
concentration phenomena, where cell density undergoes rapid growth near specific concentration points or along certain curves. Such growth
can result in singular, spiky structures and lead to finite-time blowups. 

Our investigation focuses on the dynamics of the Patlak-Keller-Segel chemotaxis system and its two-species extensions. In the latter case,
different species may exhibit distinct chemotactic sensitivities, giving rise to very different rates of cell density growth. Such a 
situation may be extremely challenging for numerical methods as they may fail to accurately capture the blowup of the slower-growing species
mainly due to excessive numerical dissipation.

In this paper, we propose a hybrid finite-difference-particle (FDP) method, in which a particle method is used to solve the chemotaxis
equation(s), while finite difference schemes are employed to solve the chemoattractant equation. Thanks to the low-dissipation nature of the
particle method, the proposed hybrid scheme is particularly adept at capturing the blowup behaviors in both one- and two-species cases. The
proposed hybrid FDP methods are tested on a series of challenging examples, and the obtained numerical results demonstrate that our hybrid
method can provide sharp resolution of the singular structures even with a relatively small number of particles. Moreover, in the
two-species case, our method adeptly captures the blowing-up solution for the component with lower chemotactic sensitivity, a feature not
observed in other works.
\end{abstract}



\section{Introduction}
We consider the two-dimensional (2-D) two-species Patlak-Keller-Segel (PKS) type chemotaxis system:
\begin{subnumcases}{\label{KS}}
(\rho_1)_t+\chi_1\nabla\cdot(\rho_1\nabla c)=\nu_1\Delta\rho_1,\label{KS1}\\
(\rho_2)_t+\chi_2\nabla\cdot(\rho_2\nabla c)=\nu_2\Delta\rho_2,\label{KS2}\\
\tau c_t=\nu\Delta c+\gamma_1\rho_1+\gamma_2\rho_2-\zeta c,\label{KS3}
\end{subnumcases}
where $\bm x=(x,y)^\top\in\Omega\subset\mathbb R^2$ are the spatial variables, $t$ is time, $\rho_1$ and $\rho_2$ represent the cell
densities of two non-competing species, $c$ is the concentration of the chemoattractant, $\chi_2>\chi_1>0$ are the chemotactic sensitivity
parameters, $\nu_1$, $\nu_2$, and $\nu$ are the diffusion coefficients, $\gamma_1$ and $\gamma_2$ represent chemoattractant production
rates, and $\zeta$ represents chemoattractant decay rate ($\nu_1$, $\nu_2$, $\nu$, $\gamma_1$, $\gamma_2$, and $\zeta$ are positive
constants). The parameter $\tau$ is either 0 or 1, corresponding to the case of parabolic-elliptic or parabolic-parabolic coupling,
respectively. The system \eref{KS} is considered subject to certain initial data and the homogeneous Neumann boundary conditions.

The original single-species PKS chemotaxis model, which is obtained from the system \eref{KS} by setting $\rho_2\equiv0$, was first
introduced in \cite{Keller1970,KS711,Pat} and later extended to the two-species system \eref{KS} in \cite{Wolansky2002}. Analytical
investigations followed in \cite{Herrero1996,Herrero1997,Othmer1997,Senba2001,Hillen2009,Conca2011,Espejo2012,Espejo2013}. Depending on the
initial data, boundary conditions, and parameter values, the solution may either converge to a constant steady state or develop singular
structures (blow up in finite time). In blowup scenarios, different blowup time scales are anticipated for the two densities, with $\rho_2$
possibly blowing up faster than $\rho_1$. However, as proved in \cite{Espejo2012,Espejo2013}, simultaneous blowup is the only possibility in
the parabolic-elliptic case $(\tau=0)$, even though different scalings may exist for the two species. This poses a significant challenge in
developing accurate and robust numerical methods for \eref{KS}.

To the best of our knowledge, the first numerical method for \eref{KS} was introduced in \cite{Kurganov2014}, where the system \eref{KS} was
solved a second-order hybrid finite-volume-finite-difference method proposed in \cite{Chertock2018} for the single-species PKS system. It
was observed in \cite{Kurganov2014} that when $\chi_2\gg\chi_1$, the cell densities $\rho_1$ and $\rho_2$ exhibit different blowup
behaviors: $\rho_2$ develops a $\delta$-type singularity, whereas $\rho_1$ undergoes algebraic blowup. This observation was numerically
validated by conducting a very careful mesh refinement study. However, the obtained results were misleading since only algebraic
growth---not blowup---of $\rho_1$ could be observed due to inadequate resolution (even on a very fine mesh). A fourth-order hybrid
finite-volume-finite-difference method developed in \cite{Chertock2018} helped only slightly to enhance the resolution. More recently, in
\cite{Chertock2019} an adaptive moving mesh (AMM) finite-volume semi-discrete upwind method was introduced to enhance the approximation of
singular structures of $\rho_1$ in the two-species chemotaxis system. Despite this advancement, the numerical experiments reported in
\cite{Chertock2019} show that even the AMM method failed to accurately capture the $\rho_1$-component of the blowing-up solution when
$\chi_2\gg\chi_1$.

In this paper, we offer a hybrid finite-difference-particle (FDP) approach for capturing the singular behavior of solutions of \eref{KS}. In
order to achieve a high resolution of both components of the blowing up solution ($\rho_1$ and $\rho_2$), we employ a sticky particle method
introduced in \cite{BJ,Chertock2007} in the context of pressureless gas dynamics equations, to the chemotaxis equations \eref{KS1} and
\eref{KS2}, while applying a second-order finite-difference (FD) semi-discrete scheme for solving the chemoattractant equation \eref{KS3}.

Sticky particle methods belong to a class of deterministic particle methods, which provide a diffusion-free (or low diffusion) alternative
to Eulerian methods for a variety of time-dependent PDEs; see, e.g., \cite{Chertock2017a,Chertock2007} and references therein. In these
methods, the numerical solution is sought as a linear combination of Dirac distributions ($\delta$-functions), whose positions and
coefficients represent the locations and weights of the particles, respectively. The solution is then found by following the time evolution 
of the locations and weights of the particles according to a system of ODEs obtained by considering a weak formulation of the studied PDEs.
Even though over the years deterministic particle methods were mainly used to numerically solve purely convective equations (see, e.g.,
\cite{CotKou,Puc,Rav} and references therein), the range of applicability of these methods was extended to treat convection-diffusion and
other types of equations (see, e.g., \cite{BCK05,CDKK,CLdisp,CotKou,CKM13,Degond1989,DMus,Eldredge2002,MGP,Puc}). Several approaches have
been explored for treating diffusion terms in particle methods. One of the widely used treatments is the random walk approach
\cite{Cho,Griebel00}, in which the diffusion operator is naturally approximated by adding a Wiener process to the motion of each particle.
This way, the diffusion only affects the positions of particles---not their weights. A drawback of the random walk method, like other
stochastic methods, is that its accuracy is very low, and reasonable resolutions can be achieved only if a very large number of
realizations is used, which would make the particle method for the PKS-type systems extremely inefficient. Among the deterministic
approaches, the one we have successfully implemented in this paper is a weighted particles method \cite{Degond1989,DMus,MGP}, which
is both efficient and accurate. In this method, the diffusion is first approximated by an integral operator, which is then treated as a
source term that affects the particle weights---not their positions.

As mesh-free, particle methods are quite flexible as the particle positions are self-adapted to the local flow. This, however, comes at the
expense of the regularity of the particle distribution: inter-particle distances typically change in time, and just as particles may cluster,
they may also spread away from each other. This may lead to a poor resolution of the computed solution and a low
efficiency of the particle method. The latter is related to the fact that the time step for the ODE solver used to evolve the particle
system in time depends generally on the distance between the neighboring particles. Thus, the success of deterministic particle methods
relies upon an accurate and efficient particle redistribution algorithm, which ensures that different parts of the computed solution are
adequately resolved. A large variety of re-meshing techniques has been proposed (see, e.g.,
\cite{Chertock2017a,CDTM,CKpart,CKL_2CH,CLdisp,CLP13,CKM13,Eldredge2002,Pes} to name a few) including particle merger approach used in
\cite{CKL_2CH,CLP13,CKM13}. Particle merger becomes crucial in cases when the solution naturally develops sharp spikes or even
$\delta$-type singularities. In these cases, the clustering particles must be merged into heavier particles located in the center of mass of
the merged particles: this is the key feature of the sticky particle methods, which were successfully applied to the pressureless gas 
dynamics equations \cite{Chertock2007}, dusty gas flows \cite{Chertock2017}, and the PKS-type system in this paper.

The success of the proposed hybrid FDP method hinges on the accurate and efficient computations of the values of $c$ at the particle 
locations and the projection of the particle approximations onto the nodes of the FD grid. The former is achieved with the help of a global 
piecewise linear approximation of $c$ (this technique is borrowed from the finite-volume methods; also see the finite-volume-particles 
methods in \cite{Chertock2017,Chertock2007}), while for the latter one, we introduce a new particle-grid projection approach.

The paper is organized as follows. In \S\ref{sec2}, we give a brief overview of the weighted particle method for convection-diffusion
equations. In \S\ref{sec3}, we introduce a hybrid FDP method, whose components---a novel sticky particle method for the cell densities
equations and a second-order FD scheme for the chemoattractant equations---are described in \S\ref{sec31} and \S\ref{sec32}, respectively.
Details of the projection between the particle and grid data, which needs to be performed in the hybridization algorithm, are discussed in
\S\ref{sec33}. In \S\ref{sec4}, we demonstrate the performance of the hybrid FDP method on a number of challenging numerical examples
designed to demonstrate the capability of the proposed method to resolve the blowup solution behavior with high resolution. Finally, we
conclude the paper in \S\ref{sec5} by summarizing our contributions.

\section{Weighted Particle Method: an Overview}\label{sec2}
In this section, we briefly describe the weighted particle method for 2-D convection-diffusion equations. We consider the following model
problem:
\begin{equation}
\rho_t+\nabla\cdot(\rho\Bu)=\nu\Delta\rho,\quad\bm x\in\Omega\subset\mathbb R^2,~~t>0,
\label{2.1}
\end{equation}
subject to the initial data
\begin{equation}
\rho(\Bx,0)=\rho_0(\Bx).
\label{2.2}
\end{equation}
Here, $\rho(\bm x,t)$ is an unknown function, $\Bu(\bm x,t)=(u(\bm x,t),v(\bm x,t))^\top$ is the given velocity field, $\nu$ is a positive 
diffusion coefficient, and $\rho_0(\bm x)$ is a given function. Notice that equations \eref{KS1} and \eref{KS2} would read as the 
convection-diffusion equation \eref{2.1} if one replaces $\nabla c$ there with a given velocity field. Dependence of \eref{KS1} and 
\eref{KS2} on the chemoattractant equation \eref{KS3} makes the development of particle methods for \eref{KS1} and \eref{KS2} substantially 
more challenging, as discussed in \S\ref{sec31}, where we develop a new sticky particle method for chemotaxis equations.

We aim to find an approximate solution of \eref{2.1}--\eref{2.2} as a linear combination of Dirac $\delta$-functions,
\begin{equation}
\widehat\rho\,(\Bx,t)=\sum_{i=1}^Nw_i(t)\,\delta(\Bx-\Bx_i(t)),
\label{2.3}
\end{equation}
where $\Bx_i(t)$ represents the location of the $i$-th particle, $w_i(t)$ denotes its weight, and $N$ is the total number of particles. We
would like to emphasize that this concept of particles represents a mathematical abstraction, differing from the physical concept of
particles of specific materials.

We begin with the initialization of the particle approximation \eref{2.3}. To this end, we divide $\Omega$ into a set of non-overlapping
subdomains $\{\Omega_i(0)\}_{i=1}^N$:
\begin{equation*}
\Omega=\bigcup_{i=1}^N\Omega_i(0)\quad\mbox{and}\quad\Omega_i(0)\cap\Omega_j(0)=\emptyset~~\forall i\ne j.
\end{equation*}
We set the initial location of the $i$-th particle, $\bm x_i(0)$, to be the center of mass of $\Omega_i(0)$. The corresponding initial
weight is then given by $\int_{\Omega_i(0)}\rho(\Bx,0)\,{\rm d}\Bx$, which we approximate using the midpoint rule resulting in
\begin{equation}
w_i(0)=\rho(\Bx_i(0),0)\,|\Omega_i(0)|.
\label{2.4}
\end{equation}

Following the lines of \cite{Chertock2017a} (also see \cite{Degond1989,DMus,MGP}), we replace the Laplacian operator $\Delta\rho$ on the 
right-hand side of \eref{2.1} with its integral approximation so that we obtain the equation
\begin{equation}
\rho_t(\bm x,t)+\nabla\cdot\big(\rho(\bm x,t)\Bu(\bm x,t)\big)=
\frac{\nu}{\sigma^2}\iint\limits_{\mathbb R^2}\eta_\sigma(\bm x-\bm y)\big(\rho(\bm y,t)-\rho(\bm x,t)\big)\,{\rm d}\bm y,
\label{2.5}
\end{equation}
where $\sigma$ is a small positive number and
\begin{equation}
\eta_\sigma(\bm x):=\frac{1}{\sigma^2}\eta\Big(\frac{\bm x}{\sigma}\Big),\quad\eta(\bm x)=\frac{4}{\pi}e^{-|\bm x|^2},
\label{2.6}
\end{equation}
and substitute \eref{2.3} into a weak form of \eref{2.5}--\eref{2.6} to end up with the following system of ODEs for the particle locations
$\Bx_i(t)$, weights $w_i(t)$, and subdomain sizes $|\Omega_i(t)|$:
\begin{equation}
\begin{dcases} 
\frac{{\rm d}}{{\rm d}t}\,\Bx_i(t)=\Bu(\bm x_i(t),t)=:\Bu_i(t)=(u_i(t),v_i(t))^\top,\\ 
\frac{{\rm d}}{{\rm d}t}\,w_i(t)=\frac{\nu}{\sigma^2}\sum_{j=1}^N\eta_\sigma\big(\Bx_i(t)-\Bx_j(t)\big)
\big[w_j(t)|\Omega_i(t)|-w_i(t)|\Omega_j(t)|\big]=:\beta_i(t),\\
\frac{{\rm d}}{{\rm d}t}\,|\Omega_i(t)|=\nabla\cdot\Bu(\Bx_i(t),t)\,|\Omega_i(t)|.
\end{dcases}
\label{2.7}
\end{equation}
The initial conditions for this ODE system are specified in \eref{2.4}.

The ODE system \eref{2.7} is to be numerically integrated by an appropriate ODE solver. The time step $\dt$ should depend on the distance
between nearby particles as the stability condition imposes that no particle trajectories should intersect within the time interval
$[t,t+\dt]$.  

In order to quantify this requirement, let us consider $i$-th and $j$-th particles for some $i$ and $j$ and assume, for the sake of 
simplicity,  that the system \eref{2.7} is numerically solved by the first-order forward Euler method. In this case, the particle
trajectories are straight lines and can be described by the following parametric form with parameters $\tau_i$ and $\tau_j$:
$$
\left\{\begin{aligned}&x=x_i(t)+\tau_iu_i(t),\\&y=y_i(t)+\tau_iv_i(t),\end{aligned}\right.\quad\mbox{and}\quad
\left\{\begin{aligned}&x=x_j(t)+\tau_ju_j(t),\\&y=y_j(t)+\tau_jv_j(t).\end{aligned}\right.
$$
It is easy to verify that these two straight lines intersect at the point that corresponds to
\begin{equation}
\begin{aligned}
&\tau_i=\frac{v_i(t)(x_j(t)-x_i(t))-u_i(t)(y_j(t)-y_i(t))}{u_i(t)v_j(t)-u_j(t)v_i(t)},\\
&\tau_j=\frac{v_j(t)(x_j(t)-x_i(t))-u_j(t)(y_j(t)-y_i(t))}{u_i(t)v_j(t)-u_j(t)v_i(t)},
\end{aligned}
\label{2.8}
\end{equation}
unless the denominator $u_i(t)v_j(t)-u_j(t)v_i(t)$ vanishes, in which case these two lines are parallel. Notice that the $i$-th particle
will arrive in this intersection point at time $t+\tau_i$ provided $\tau_i>0$, while the $j$-th particle will be there at time $t+\tau_j$
provided $\tau_j>0$. This leads to the following time step constraint:
\begin{equation}
\dt\le\left\{\begin{aligned}
&\min_{(i,j):\tau_i>0,\,\tau_j>0}\left\{\max\left(\tau_i,\tau_j\right)\right\}&&\mbox{if}~\exists\,(i,j): \tau_i>0,\,\tau_j>0,\\
&\infty,&&\mbox{otherwise}.
\end{aligned}\right.
\label{2.9}
\end{equation}
\begin{rmk}
In practice, in order to prevent division by very small numbers, one has to use the formulae in \eref{2.8} only when
$|u_i(t)v_j(t)-u_j(t)v_i(t)|>\mu$, where $\mu$ is a small positive number.
\end{rmk}

Additional time step constraints are obtained by requiring the weights $w_i$ to remain nonnegative for all $i$ (this is extremely important
when $\rho$ must remain nonnegative as in the studied case of chemotaxis systems), the size of subdomains $\Omega_i$ not to decay too fast
(for instance, one may restrict their decay within one-time step by a factor of 2), and no particle to propagate too far in one time step.
In the case of forward Euler time discretization, these restrictions read as
\begin{equation}
\dt\le\left\{\begin{aligned}
&\min_{i:\beta_i(t)<0}\left\{-\frac{w_i(t)}{\beta_i(t)}\right\}&&\mbox{if}~\max_i\left\{-\frac{w_i(t)}{\beta_i(t)}\right\}>0,\\
&\infty,&&\mbox{otherwise},
\end{aligned}\right.
\label{2.10}
\end{equation}
\begin{equation}
\dt\le\left\{\begin{aligned}
&\min_{i:\nabla\cdot\Bu(\Bx_i(t),t)<0}\left\{-\frac{1}{2\,\nabla\cdot\Bu(\Bx_i(t),t)}\right\}&&\mbox{if}~\max_i\left\{-\frac{1}{2\,\nabla\cdot\Bu(\Bx_i(t),t)}\right\}>0,\\
&\infty,&&\mbox{otherwise},
\end{aligned}\right.
\label{2.11}
\end{equation}
and
\begin{equation}
\dt\le\frac{\sqrt{\min_i|\Omega_i(0)|}}{\max_i\left(\max\{|u_i|,|v_i|\}\right)},
\label{2.12}
\end{equation}
respectively.
\begin{rmk}
The time restrictions \eref{2.9}--\eref{2.11} are still valid if a high-order strong stability preserving (SSP) Runge-Kutta or multistep
explicit method \cite{GKS,Gottlieb2001} is used instead of the forward Euler method for time integration, since the SSP methods can be
expressed as convex combinations of forward Euler steps.
\end{rmk}

\section{Hybrid FDP Method for \eref{KS}}\label{sec3}
We now turn our attention to the PKS-type system \eref{KS} and introduce an FDP method for its numerical solution. Specifically, the
chemotaxis equations \eref{KS1} and \eref{KS2} are discretized using a sticky particle method (\S\ref{sec31}), while the chemoattractant
concentration equation \eref{KS3} is approximated using a FD scheme (\S\ref{sec32}). The coupling between the two ingredients of the FDP
methods is carried out using the projection techniques (\S\ref{sec33})

\subsection{Sticky Particle Method for Equations \eref{KS1} and \eref{KS2}}\label{sec31}
Recall that the chemotaxis equations \eref{KS1} and \eref{KS2} take the form of the convection-diffusion equation \eref{2.1} with $\rho$
being either $\rho_1$ or $\rho_2$, $\nu=\nu_1$ or $\nu=\nu_2$, and $\bm u=\chi_1\nabla c$ or $\bm u=\chi_2\nabla c$. The
corresponding particle solutions for $\rho_1$ and $\rho_2$ are sought in the form
\begin{equation}
\widehat\rho_1(\Bx,t)=\sum_{i=1}^{N_1}w^1_i(t)\,\delta(\Bx-\Bx^1_i(t))\quad\mbox{and}\quad
\widehat\rho_2(\Bx,t)=\sum_{i=1}^{N_2}w^2_i(t)\,\delta(\Bx-\Bx^2_i(t))
\label{3.1}
\end{equation}
with $\Bx^1_i(t)$ and $\Bx^2_i(t)$ being each species particle locations, $w^1_i(t)$ and $w^2_i(t)$ their corresponding weights, and $N_1$
and $N_2$ the total number of each species particles. We denote by $\Omega^1_i(t)$ and $\Omega^2_i(t)$ the subdomains occupied by the
corresponding particles.

In principle, one can apply the weighted particle method from \S\ref{sec2} to \eref{KS1} and \eref{KS2}, in which the particle locations,
their weights, and sizes of the corresponding subdomains evolve in time according to \eref{2.7}. However, when the solutions of \eref{KS1}
and \eref{KS2} start developing spiky structures, the particles start clustering in the regions of large chemoattractant gradient and
hence, using a constant $\sigma$ in \eref{2.7} leads to an inaccurate approximation of the diffusion operator. One therefore has to use
variable values of $\sigma$, which depend on the distance between the particles. In particular, we use
$\sigma^k_{ij}=\sqrt{\big(|\Omega^k_i|+|\Omega^k_j|\big)/2}$, $k=1,2$, so that the system \eref{2.7} for $\rho_1$ and $\rho_2$ reads as
\begin{equation}
\begin{dcases} 
\frac{{\rm d}}{{\rm d}t}\,\Bx^k_i=(u^k_i,v^k_i)^\top,\quad k=1,2,\\ 
\frac{{\rm d}}{{\rm d}t}\,w^k_i=\nu\sum_{j=1}^{N_k}\frac{1}{(\sigma^k_{ij})^2}\,\eta_{\sigma^k_{ij}}\big(\Bx^k_i-\Bx^k_j\big)
\big[w^k_j|\Omega^k_i|-w^k_i|\Omega^k_j|\big],\quad k=1,2,\\
\frac{{\rm d}}{{\rm d}t}\,|\Omega^k_i|=r^k_i|\Omega^k_i|,\quad k=1,2.
\end{dcases}
\label{3.2}
\end{equation}
Here, $(u^k_i,v^k_i)^\top\approx\chi_k\nabla c(x^k_i(t),y^k_i(t),t)$ and $r^k_i\approx\chi_k\Delta c(x^k_i(t),y^k_i(t),t)$, which will be
obtained in \S\ref{sec331} below. Notice that all of the indexed quantities in \eref{3.2} depend on $t$, but we have omitted this dependence
for the sake of brevity. Time dependence of indexed quantities will be also omitted below unless it is necessary to emphasize it in the
discussion context. 

It should be observed that choosing variable $\sigma^k_{ij}$ is not sufficient to make numerical integration of the ODE system \eref{3.2}
practically feasible as the time step restriction \eref{2.9} would lead to $\dt\to0$, when $\rho_1$ and $\rho_2$ approaching the point mass
concentration.

In order to overcome this difficulty, we present a new sticky particle method for the chemotaxis equations \eref{KS1} and \eref{KS2}. The
basic idea of sticky particle methods is to coalesce clustering particles into a ``heavier'' particle located at their center of mass. In
this paper, we intoduce a new particle merger strategy, which is implemented in two steps. Namely, at each time level $t$, we first identify
pairs of particles whose trajectories are about to intersect by the next time level $t+\dt$ and merge them before the time evolution step.
We then evolve the particles from $t$ to $t+\dt$ and coalesce those particles, which have clustered upon the completion of the evolution
step (the second merger step was used in \cite{Chertock2017,CKL_2CH,Chertock2007,CLP13,CKM13}).

In order to implement the merger, we introduce an auxiliary ``merger'' Cartesian grid consisting of small cells of the size about
$\frac{1}{4}\min_i|\Omega_i(0)|$ and assume that at a certain time level $t$ no ``merger'' cell contain more than one particle. We then
compute the time step $\dt$ based on the two-species version of the time step restrictions \eref{2.10}--\eref{2.12}, but not \eref{2.9},
which, as we have mentioned, may be impractically too restrictive. Equipped with this $\dt$, we proceed with the following steps.

\paragraph{Merger Step 1.} We identify those pairs of the same species particles, whose trajectories are about to intersect before time
$t+\dt$. Let us assume that $\{(\bm x_i^k,w_i^k,|\Omega_i^k|),(\bm x_j^k,w_j^k,|\Omega_j^k|)\}$ is one of such pairs ($k=1$ or 2), that is,
$\tau_i^k>0$, $\tau_j^k>0$, and $\max\{\tau_i^k,\tau_j^k\}<\dt$, where $\tau_i^k$ and $\tau_j^k$ are given by the two-species version of
\eref{2.8}. We then add the weights of these two particles, combine their subdomains, and replace the two particles with one heavier
particle located at the center of mass of the replaced particles:
$$
\left\{(\bm x_i^k,w_i^k,|\Omega_i^k|),(\bm x_j^k,w_j^k,|\Omega_j^k|)\right\}\quad\longrightarrow\quad
(\widetilde{\bm x}^k,\widetilde w^k,|\widetilde\Omega^k|)
$$
with
$$
\widetilde{\bm x}^k=\frac{w_i^k\bm x_i^k+w_j^k\bm x_j^k}{w_i^k+w_j^k},\quad\widetilde w^k=w_i^k+w_j^k,\quad
|\widetilde\Omega^k|=|\Omega_i^k|+|\Omega_j^k|.
$$
Once the $i$-th and $j$-th particles of species $k$ have been merged, the total number of particles for this species reduces by one, the
remaining particles are re-numbered, and Merger Step 1 is repeated for each species to ensure that no particle trajectories intersect within
the current time step.
\begin{rmk}
Notice that the process of searching for pairs of particles with potentially intersecting trajectories is computationally expensive if 
performed in a straightforward manner---looping through all particle pairs of species $k$ is ${\cal O}(N_k^2)$ expensive. However, the
search mechanism can be made efficient by introducing another auxiliary Cartesian grid with cells of size $\min_i|\Omega_i(0)|$ and using
the fact that the time step restriction \eref{2.12} ensures that particles cannot propagate too far within one time step. One can then check
the particles located in nearby auxiliary cells. This alternative implementation is only ${\cal O}(N_k)$ expensive for $k=1,2$.
\end{rmk}

\paragraph{Time Evolution.} The particle solution is evolved from time $t$ to $t+\dt$ by numerically solving \eref{3.2}. Recall that the
system \eref{KS} is solved subject to the homogeneous Neumann boundary conditions, which imply that no particles should leave the
computational domain $\Omega$. This is, however, not automatically guaranteed and we therefore develop a ``pull-back'' strategy of
relocating the outside particles back on the domain boundary $\partial\Omega$ as follows: we move the outside particle to the closest point
on $\partial\Omega$. Once all of the evolved particles are located in $\Omega$, one has to coalesce the clustering particles according to
Merger Step 2.

If a Runge-Kutta ODE solver is employed (as in all of the numerical experiments presented in \S\ref{sec4}), then both ``pull-back'' and
merger algorithms should be implemented at the end of each Runge-Kutta stage.

\paragraph{Merger Step 2.} We find all of the auxiliary ``merger'' cells containing more than one particle from the same species and then
merge them according to the following procedure. Let $C_{\rm mer}$ be one of such cells with two or more particles of species $k$ ($k=1$ or
2). Then, the particles of species $k$ located in $C_{\rm mer}$ are merged into a new particle,
$(\widetilde{\bm x}^k,\widetilde w^k,|\widetilde\Omega^k|)$, located at the center of mass of the replaced particles and their weights and
subdomain sizes are summed up:
\begin{equation*}
\widetilde{\bm x}^k=\frac{\sum\limits_{i: \Bx_i^k\in C_{\rm mer}}\hspace*{-0.4cm}w_i^k\Bx_i^k}
{\sum\limits_{i: \Bx_i^k\in C_{\rm mer}}\hspace*{-0.4cm}w_i^k},\qquad
\widetilde w^k=\hspace*{-0.4cm}\sum\limits_{i: \Bx_i^k\in C_{\rm mer}}\hspace*{-0.4cm}w_i^k,\qquad|
\widetilde\Omega^k|=\hspace*{-0.4cm}\sum\limits_{i: \Bx_i\in C_{\rm mer}}\hspace*{-0.4cm}|\Omega_i^k|.
\end{equation*}
\begin{rmk}
We would like to point out that, by construction, the sticky particle method preserves the positivity of $\rho_1$ and $\rho_2$, which is a
crucial property for the stability of numerical methods for PKS-type systems, as it was show in \cite{Chertock2008}.
\end{rmk}
   
\subsection{Finite-Difference Scheme for Equation \eref{KS3}}\label{sec32}
In this section, we describe the second ingredient of our hybrid method---the FD scheme, which is used to numerically solve the
chemoattractant concentration equation \eref{KS3}. We restrict our consideration to the case of a rectangular computational domain $\Omega$.
For general domains, one can still use FD schemes, but treatment of boundary conditions and grid points near the boundary becomes quite
delicate (this is outside the scope of current paper).

We first split $\Omega$ into uniform rectangular cells $C_{\ell,m}$ of dimensions $\dx$ and $\dy$ and denote the cell centers by
$\bm x_{\ell,m}=(x_{\ell,m},y_{\ell,m})$. A second-order FD discretization of \eref{KS3} (or, a semi-discretization in the parabolic case
with $\tau=1$) reads as
\begin{equation}
\tau\frac{{\rm d}}{{\rm d}t}c_{\ell,m}=\nu\bigg[\frac{c_{\ell+1,m}-2c_{\ell,m}+c_{\ell-1,m}}{(\dx)^2}+
\frac{c_{\ell,m+1}-2c_{\ell,m}+c_{\ell,m-1}}{(\dy)^2}\bigg]+\gamma_1(\rho_1)_{\ell,m}+\gamma_2(\rho_2)_{\ell,m}-\zeta c_{\ell,m},
\label{3.3}
\end{equation}
where $c_{\ell,m}:\approx c(\Bx_{\ell,m},t)$.

If $\tau=1$, then one need to numerically integrate the extended ODE system \eref{3.2}--\eref{3.3}. Otherwise, if $\tau=0$, \eref{3.3} is a
linear algebraic system with respect to $\{c_{\ell,m}\}$, which needs to be solved upon completion of each Runge-Kutta stage.

\subsection{Projection Between the Particle and Grid Data}\label{sec33}
It should be obsevred that one has to project the data from/to the particle locations $\Bx^k_i$ to/from the grid nodes $\Bx_{\ell,m}$, when
a meshless particle method is combined with a grid-based FD scheme. We introduce these mapping procedures in \S\ref{sec331} and
\S\ref{sec332}. 

\subsubsection{Computation of $\nabla c$ and $\Delta c$ at Particle Locations}\label{sec331}
In order to evolve particles in time according to \eref{3.2}, one has to compute velocities $(u_i^k,v_i^k)^\top$ and their divergences
$r_i^k$ at the particle locations. Since these quantities are obtained using $\nabla c$ and $\Delta c$ and since the values $c_{\ell,m}$ are
computed at the grid points, a projection of the grid data on the particle location has to be carried out. This can be done as follows.

First, we compute the point values of $c_x$, $c_y$, and $\Delta c$ at the grid points using the second-order FDs:
\begin{equation*}
\begin{aligned}
&(c_x)_{\ell,m}=\frac{c_{\ell+1,m}-c_{\ell-1,m}}{2\dx},\quad(c_y)_{\ell,m}=\frac{c_{\ell,m+1}-c_{\ell,m-1}}{2\dy},\\
&(\Delta c)_{\ell,m}=\frac{c_{\ell+1,m}-2c_{\ell,m}+c_{\ell-1,m}}{(\dx)^2}+\frac{c_{\ell,m+1}-2c_{\ell,m}+c_{\ell,m-1}}{(\dy)^2},
\end{aligned}
\end{equation*}
then construct global (in space) piecewise linear interpolants for $c_x$, $c_y$, and $\Delta c$, which is in every cell $C_{\ell,m}$ defined
by
$$
\begin{aligned}
&\widetilde{c_x}(x,y)=(c_x)_{\ell,m}+\frac{(c_x)_{\ell+1,m}-(c_x)_{\ell-1,m}}{2\dx}(x-x_{\ell,m})+
\frac{(c_x)_{\ell,m+1}-(c_x)_{\ell,m-1}}{2\dy}(y-y_{\ell,m}),\\
&\widetilde{c_y}(x,y)=(c_y)_{\ell,m}+\frac{(c_y)_{\ell+1,m}-(c_y)_{\ell-1,m}}{2\dx}(x-x_{\ell,m})+
\frac{(c_y)_{\ell,m+1}-(c_y)_{\ell,m-1}}{2\dy}(y-y_{\ell,m}),\\
&\widetilde{\Delta c}(x,y)=(\Delta c)_{\ell,m}+\frac{(\Delta c)_{\ell+1,m}-(\Delta c)_{\ell-1,m}}{2\dx}(x-x_{\ell,m})+
\frac{(\Delta c)_{\ell,m+1}-(\Delta c)_{\ell,m-1}}{2\dy}(y-y_{\ell,m}),
\end{aligned}
$$
and finally, for each particle location $\bm x_i^k$, $k=1,2$, we obtain
\begin{equation}
u^k_i=\chi_k\widetilde{c_x}(x^k_i,y^k_i),\quad v^k_i=\chi_k\widetilde{c_y}(x^k_i,y^k_i),\quad
r^k_i=\chi_k\widetilde{\Delta c}(x^k_i,y^k_i).
\label{3.4}
\end{equation}

\subsubsection{Computation of $\rho_1$ and $\rho_2$ at Grid Points}\label{sec332}
In order to evolve $c$ in time according to \eref{3.3}, one has to recover the density grid values $(\rho_1)_{\ell,m}$ and
$(\rho_2)_{\ell,m}$ from the particle distributions \eref{3.1}. This can be conducted in the following three steps.

\smallskip
\noindent
{\bf Step 1.} We compute the point values $\rho_1$ and $\rho_2$ at the particle locations $\Bx_i^k$ by dividing the particle weight $w^k_i$
by the corresponding subdomain sizes $|\Omega^k_i|$, that is, we set $\rho^k_i=w^k_i/|\Omega^k_i|$.

\smallskip
\noindent
{\bf Step 2.} For each particle satisfying $\Bx^k_i\in C_{\ell,m}$, we compute the distance from this particle to the cell center
$\Bx_{\ell,m}$:
\begin{equation}
(d^k_i)_{\ell,m}:=|\Bx^k_i-\bm x_{\ell,m}|,
\label{3.5}
\end{equation}
and evaluate the grid values of $\rho_1$ and $\rho_2$ at the cell center $\bm x_{\ell,m}$ using a distance-based weighted averaging:
\begin{equation}
(\rho_k)_{\ell,m}^*=\frac{\sum\limits_{\bm x^k_i\in C_{\ell,m}}\rho^k_i/(d^k_i)_{\ell,m}}
{\sum\limits_{\bm x^k_i\in C_{\ell,m}}1/(d^k_i)_{\ell,m}},\qquad k=1,2.
\label{3.6}
\end{equation}
If there is no particles of species $k$ in cell $C_{\ell,m}$, we set $(\rho_k)_{\ell,m}^*=0$.

\smallskip
\noindent
{\bf Step 3.}
In the cells containing no particles, the value $(\rho_k)_{\ell,m}^*=0$ may be very inaccurate. The therefore replace the values computed in
{\bf Step 2} with: 
$$
(\rho_k)_{\ell,m}=\begin{dcases}(\rho_k)_{\ell,m}^*&\mbox{if } \exists \, i: \Bx^k_i\in C_{\ell,m},\\
\begin{aligned}
&\frac{(\rho_k)_{\ell+1,m}^*+(\rho_k)_{\ell-1,m}^*+(\rho_k)_{\ell,m+1}^*+(\rho_k)_{\ell,m-1}^*}{4+2\sqrt{2}}\\
&+\frac{(\rho_k)_{\ell+1,m+1}^*+(\rho_k)_{\ell+1,m-1}^*+(\rho_k)_{\ell-1,m+1}^*+(\rho_k)_{\ell-1,m-1}^*}{4+4\sqrt{2}},
\end{aligned}
&\mbox{otherwise}.
\end{dcases}
$$
where the latter expression is the distance-based weighted averaging over neighboring cells.
\begin{rmk}
It should be observed that there are other ways to recover grid values of the computed solution from its particle distribution when the
solution is smooth; see, e.g., \cite{Chertock2017a}. However, after the spiky structure is developed, most of these methods are based on a
certain regularization of the $\delta$-functions and hence lead to a substantial decrease of the maximum values of $\rho_1$ and $\rho_2$.
The method we have used in Step 1 seems to be the only robust option.
\end{rmk}
\begin{rmk}
Note that some of the values $d^k_i$ defined in \eref{3.5} may be very small or even zero. We therefore need to desingularize the
computation in \eref{3.6} to prevent division by small numbers. This is done by replacing \eref{3.5} with
$$
d^k_i=\max\left\{d_{\min},|\Bx^k_i-\bm x_{\ell,m}|\right\},
$$
where $d_{\min}$ is a small positive number taken to be $\min(\dx,\dy)/16$ in all of the numerical examples reported in \S\ref{sec4}.
\end{rmk}

\section{Numerical Examples}\label{sec4}
In this section, we demonstrate the performance of the proposed hybrid FDP method and its capability of capturing the blowing up solutions
of the PKS-type system \eref{KS} with high resolution.

Recall that the proposed method employs three different meshes for distinct purposes: a mesh used in the FD scheme ($C_{\ell,m}$), a mesh
used as initial particle subdomains $\Omega_i(0)$, and a mesh used for Merger Step 2 ($C_{\rm mer}$). In the following numerical examples,
we use uniform meshes with $\Delta:=\dx=\dy$ and keep the ratio between these three mesh sizes fixed. Specifically, by ``numerical results
obtained on a mesh of size $\Delta$'', we imply the combination of a mesh of size $\Delta$ for the FD scheme, a mesh of size $\Delta/4$ for
initial particle subdomains, and a mesh of size $\Delta/8$ for the merger.

For time evolution, we use the three-stage third-order SSP Runge-Kutta method \cite{GKS,Gottlieb2001}. The time step is selected using the
stability restrictions introduced in \eref{2.10}--\eref{2.12}. We remind the reader that once the final time solution has been computed, the
point values of $c$ are available, while the point values of $\rho_1$ and $\rho_2$ need to be recovered as described in {\bf Steps 1 and 2}
in \S\ref{sec332}. In order to visualize the obtained solution and conduct the experimental convergence study, one may use, for instance,
the MATLAB built-in function {\tt scatteredInterpolant}.

\subsection{Parabolic-Parabolic Case ($\tau=1$)}\label{sec41}
\subsubsection*{Example 1---Accuracy Test}
The primary objective of this example is to experimentally check the accuracy of the proposed hybrid FDP method. 

We consider the system \eref{KS} with $\tau=1$, $\nu=10$, $\chi_1=5$, $\chi_2=60$, $\gamma_1=\gamma_2=\zeta=\nu_1=\nu_2=1$, and subject to
the following initial conditions:
\begin{equation*}
\rho_1(x,y,0)=\rho_2(x,y,0)=500\,{\rm e}^{-100(x^2+y^2)},\quad c(x,y,0)\equiv1,\quad(x,y)\in\Omega=[-1,1]\times[-1,1].
\end{equation*}
We compute the numerical solutions until time $t=2\times10^{-4}$ using various resolutions with $\Delta=2/15$, 2/20, 2/25, 2/30, 2/40, 2/50,
2/60, 2/80, 2/100, and 2/120. Given the numerical solutions obtained on meshes of sizes $\Delta$, $\Delta/2$, and $\Delta/4$, we compute the
$L^1$- and $L^2$-errors and estimate the corresponding experimental convergence rates using the following Runge formulae:
$$
\mbox{$L^p$-error}\approx\frac{\|(\cdot)^\frac{\Delta}{2}-(\cdot)^\frac{\Delta}{4}\|_{L^p}^2}
{\big|\|(\cdot)^\Delta-(\cdot)^\frac{\Delta}{2}\|_{L^p}-\|(\cdot)^\frac{\Delta}{2}-(\cdot)^\frac{\Delta}{4}\|_{L^p}\big|},\quad
\mbox{rate}\approx\log_2\left(\frac{\|(\cdot)^\Delta-(\cdot)^\frac{\Delta}{2}\|_{L^p}}
{\|(\cdot)^\frac{\Delta}{2}-(\cdot)^\frac{\Delta}{4}\|_{L^p}}\right),\quad p=1,2,
$$
where $(\cdot)^\Delta$ denotes the numerical results obtained on a mesh of size $\Delta$.

The computed $L^1$- and $L^2$-errors and the corresponding experimental convergence rates for $\rho_1$, $\rho_2$, and $c$ are presented in
Table \ref{tab41}, where one can see that the second order of accuracy is achieved.
\renewcommand{\tabcolsep}{3.6pt}
\begin{table}[ht!]
\centering
\begin{tabular}{ccccccccccccc}
\toprule[2pt]
&\multicolumn{4}{c}{$\rho_1$} &\multicolumn{4}{c}{$\rho_2$}&\multicolumn{4}{c}{$c$}\\
\cmidrule(lr){2-5}\cmidrule(lr){6-9}\cmidrule(lr){10-13}
$\Delta$&$L^1$-error&rate&$L^2$-error&rate&$L^1$-error&rate&$L^2$-error&rate&$L^1$-error&rate&$L^2$-error&rate\\
\midrule[0.8pt] 
2/60&2.94e-2&1.95&1.35e-1&1.88&7.70e-2&2.20&2.94e-1&2.43&1.80e-4&1.78&4.12e-4&2.00\\
2/80&1.55e-2&2.00&7.75e-2&1.85&4.52e-2&1.64&1.96e-1&1.70&1.06e-4&1.81&3.43e-4&1.65\\
2/100&1.04e-2&1.93&4.48e-2&1.93&2.19e-2&2.43&8.40e-2&2.68&7.33e-5&1.75&1.65e-4&1.93\\
2/120&6.20e-3&2.07&2.55e-2&2.10&1.93e-2&2.11&6.63e-2&2.06&4.52e-5&1.87&1.06e-4&1.99\\
\bottomrule[2pt]
\end{tabular}
\caption{\sf Example 1: The $L^1$- and $L^2$-errors and the corresponding convergence rates.\label{tab41}}
\end{table}

\subsubsection*{Example 2---Two-Species Blowup at the Center of the Domain}
\begin{figure}[ht!]
\centerline{\includegraphics[trim=0cm 2.2cm 0cm 0.9cm,clip,width=0.9\textwidth]{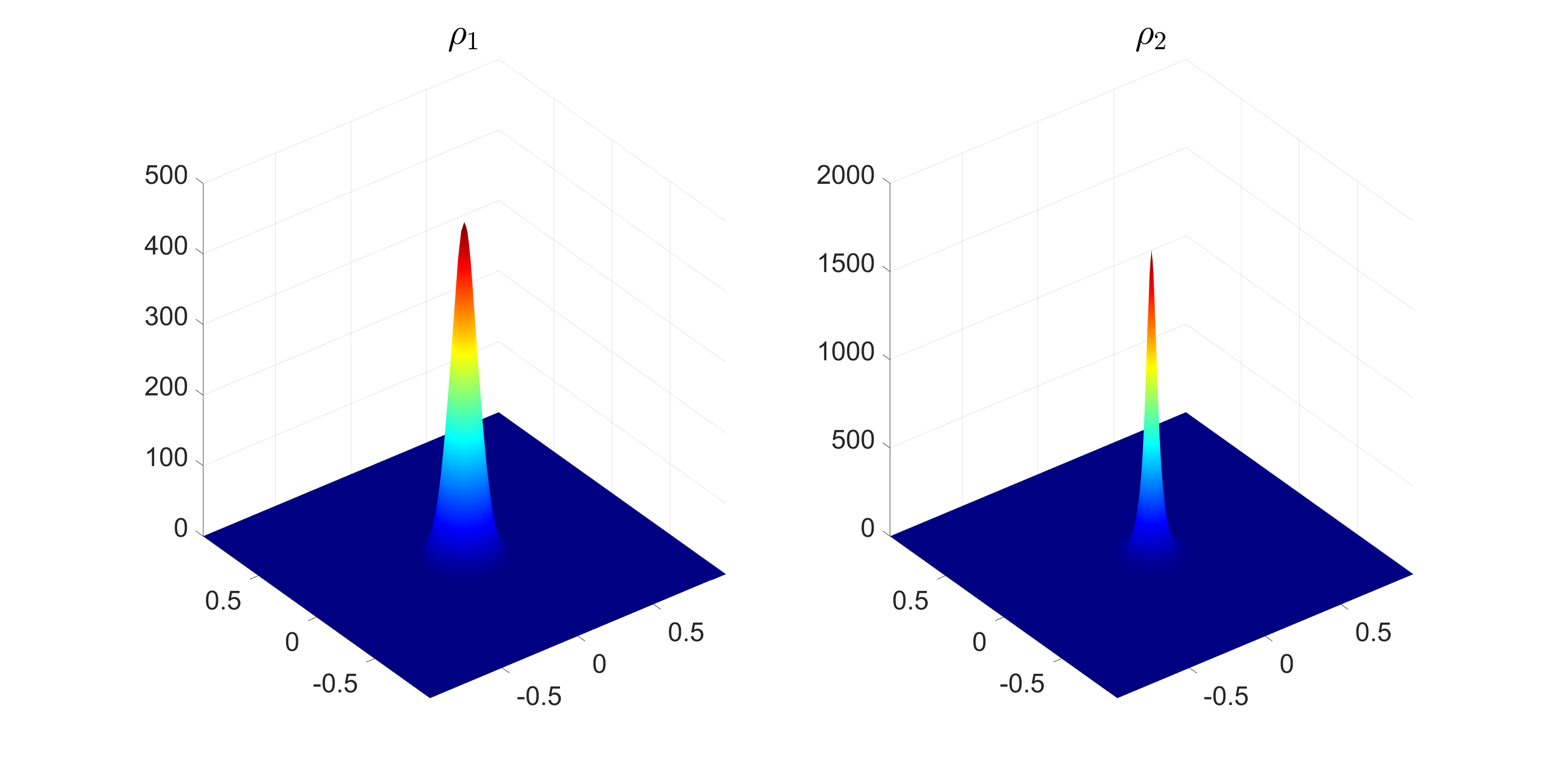}}
\vskip20pt
\centerline{\includegraphics[trim=0cm 2.2cm 0cm 0.9cm,clip,width=0.9\textwidth]{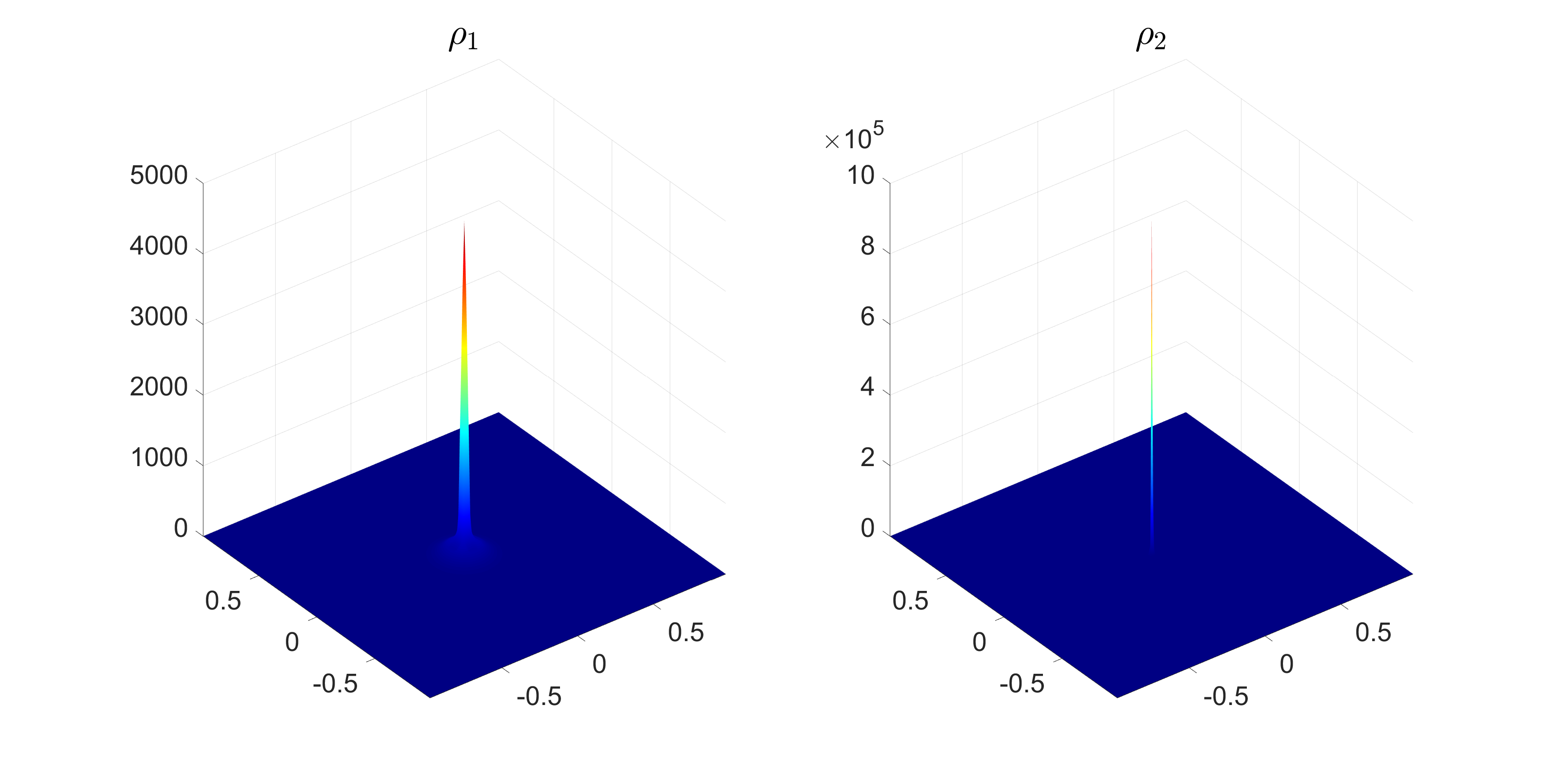}}
\caption{\sf Example 2: $\rho_1(\bm x,0.0005)$ (upper left), $\rho_2(\bm x,0.0005)$ (upper right), $\rho_1(\bm x,0.001)$ (lower left), and
$\rho_2(\bm x,0.001)$ (lower right), obtained using the hybrid FDP method with $\Delta=1/20$.\label{fig41}}
\end{figure}

In this example, we consider the same initial-boundary value problem (IBVP) as in Example 1, but compute its numerical solution until much
larger times $t=5\times10^{-4}$ and $10^{-3}$. The densities $\rho_1$ and $\rho_2$ obtained on a mesh of size $\Delta=1/20$ are plotted in
Figure \ref{fig41}. The results are comparable to the ones in \cite[Example 3]{Chertock2019} and \cite{Kurganov2014}: Both $\rho_1$ and
$\rho_2$ exhibit blowup behavior, and the solution of $\rho_2$ blows up faster.

This example is challenging in the sense that excessive numerical viscosity may smear the singularity, leading to a failure to capture the
blowup behavior. Due to the low-dissipation nature of the sticky particle method, the hybrid FDP method provides a higher resolution: The
maximum values of $\rho_1$ and $\rho_2$ are respectively $4.7355\times10^3$ and $9.6558\times10^5$, as compared to $3.5758\times10^3$ and
$2.8726\times10^5$ reported in \cite[Example 3]{Chertock2019}.

\subsubsection*{Example 3---Single-Species Blowup at the Corner of the Domain}
The third example taken from \cite{Chertock2008} is designed to demonstrate the ability of the proposed hybrid FDP method to capture the
blowup behavior away from the center of the initial Gaussian-shaped cell density. 

We consider the single-species PKS system
\begin{equation}
\left\{\begin{aligned}
&\rho_t+\nabla\cdot\left(\rho\nabla c\right)=\Delta\rho,\\
&c_t=\Delta c+\rho-c,
\end{aligned}\right.
\label{1KS}
\end{equation}
subject to following initial conditions:
\begin{equation}
\rho(x,y,0)=500\,{\rm e}^{-100[(x-0.25)^2+(y-0.25)^2]},~~c(x,y,0)\equiv0,~~(x,y)\in\Omega=[-0.5,0.5]\times[-0.5,0.5].
\label{4.2}
\end{equation}
It has been proved in \cite{Herrero1997} that when the total mass of $\rho$ is below a certain threshold, the density $\rho$ can only blow
up at the boundary of the computational domain. This is indeed the case for the IBVP \eref{1KS}--\eref{4.2}.

\begin{figure}[ht!]
\centerline{\includegraphics[trim=0cm 1.9cm 0cm 0.8cm,clip,width=0.9\textwidth]{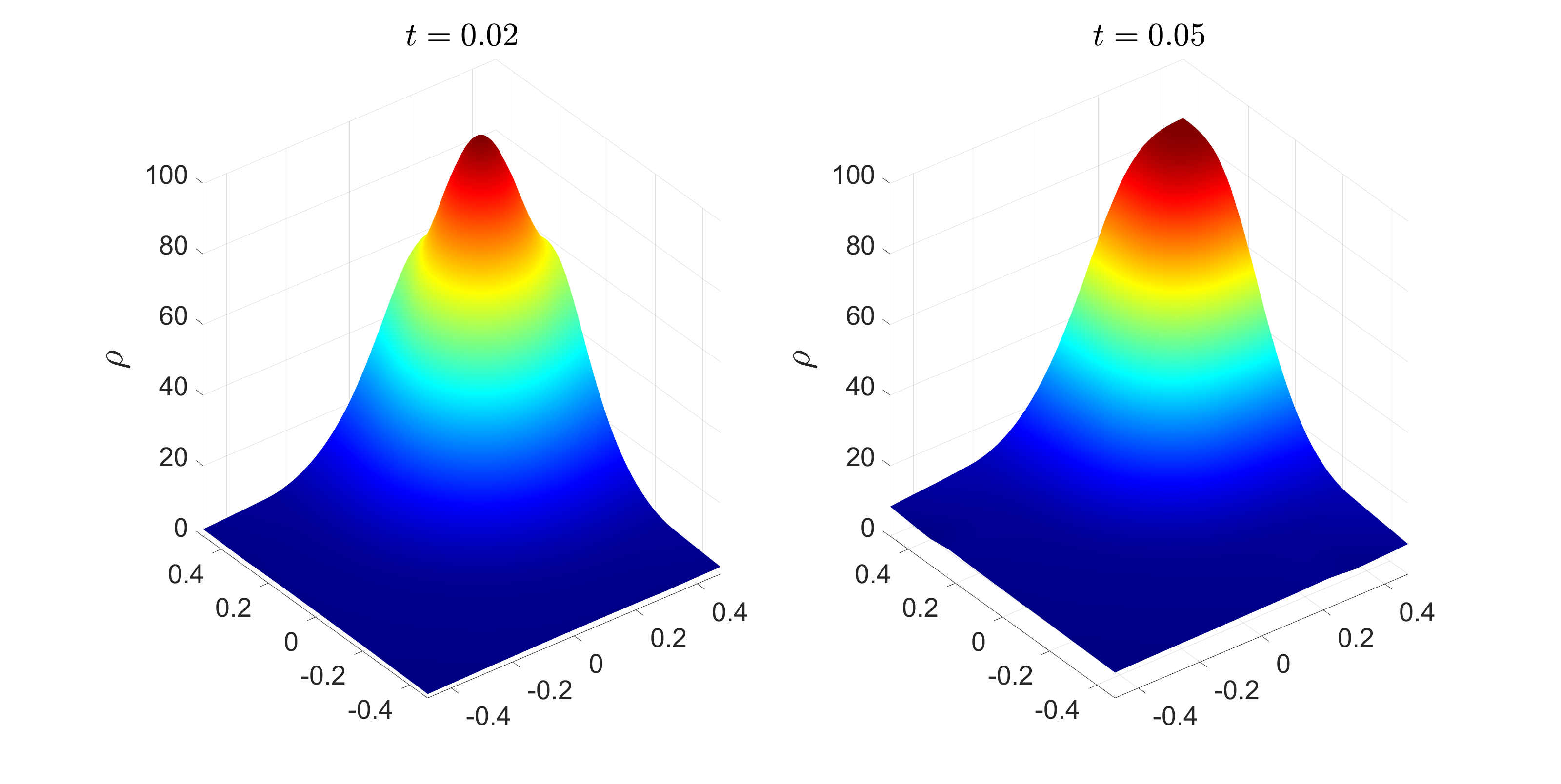}}
\vskip20pt
\centerline{\includegraphics[trim=0cm 1.9cm 0cm 0.8cm,clip,width=0.9\textwidth]{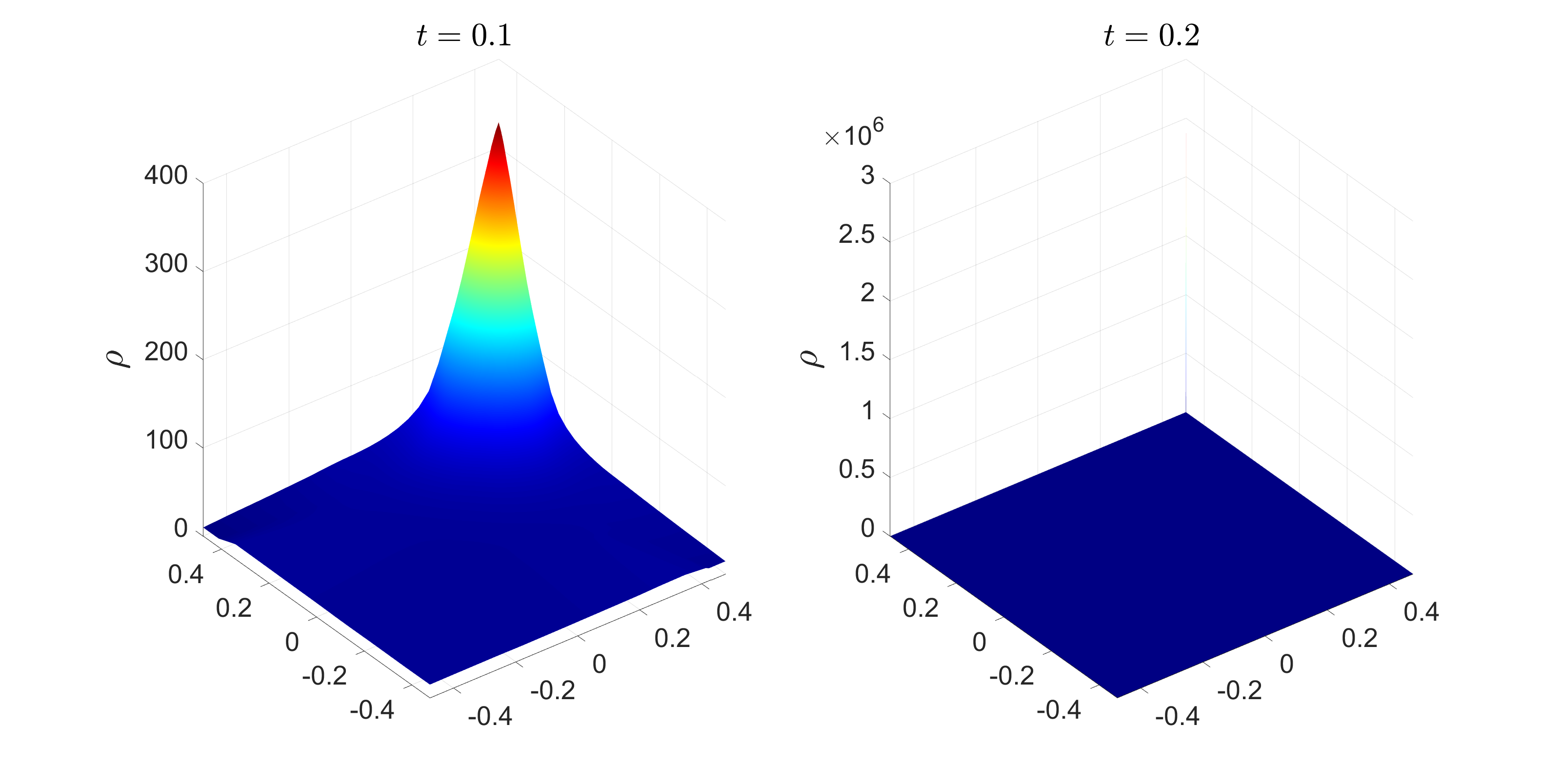}}
\caption{\sf Example 3: $\rho(\bm x,0.02)$ (upper left), $\rho(\bm x,0.05)$ (upper right), $\rho(\bm x,0.1)$ (lower left), and
$\rho(\bm x,0.2)$ (lower right), obtained using the hybrid FDP method with $\Delta=1/20$.\label{fig42}}
\end{figure}

\begin{figure}[ht!]
\centerline{\includegraphics[trim=0cm 2.3cm 0cm 2.7cm,clip,width=0.9\textwidth]{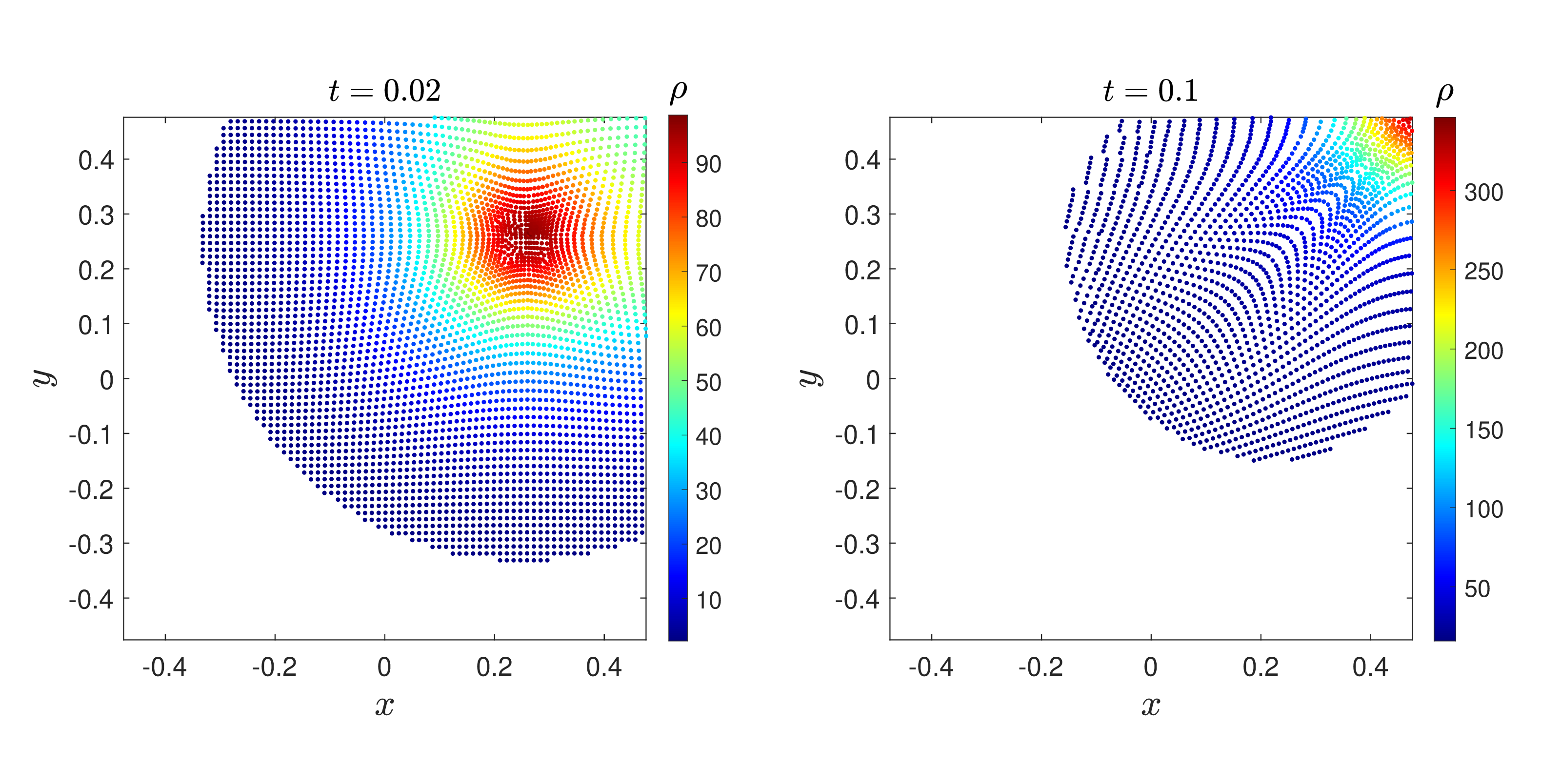}}
\caption{\sf Example 3: particle locations at $t=0.02$ (left) and $t=0.1$ (right), obtained using the hybrid FDP method with $\Delta=1/20$.
The color indicates the value of $w_i(t)/|\Omega_i(t)|$.\label{fig43}}
\end{figure}

We employ the hybrid FDP method to compute the solution on a mesh of size $\Delta=1/20$. The density $\rho$ at different times is shown in
Figure \ref{fig42}, and the particle locations $\{(x_i(t),y_i(t))\}$ at times $t=0.02$ and 0.1 are depicted in Figure \ref{fig43}. As one
can observe, the behavior of the computed solution matches the theoretical results established in \cite{Herrero1997}: The mass of $\rho$
first moves to the boundary and then concentrates at the corner where the solution blows up.

\subsubsection*{Example 4---Single-Species Blowup in the Interior of the Domain}
We consider the same IBVP as in Example 3, but with different initial conditions:
\begin{equation*}
\rho(x,y,0)=1000\,{\rm e}^{-100[(x-0.25)^2+(y-0.25)^2]},\quad c(x,y,0)\equiv0,\quad(x,y)\in\Omega=[-0.5,0.5]\times[-0.5,0.5].
\end{equation*}
In this case, the total mass of $\rho$ is greater than the threshold, and the solution is expected to blow up in the interior of the
computational domain. 

\begin{figure}[ht!]
\centerline{\includegraphics[trim=0cm 1.9cm 0cm 0.8cm,clip,width=0.9\textwidth]{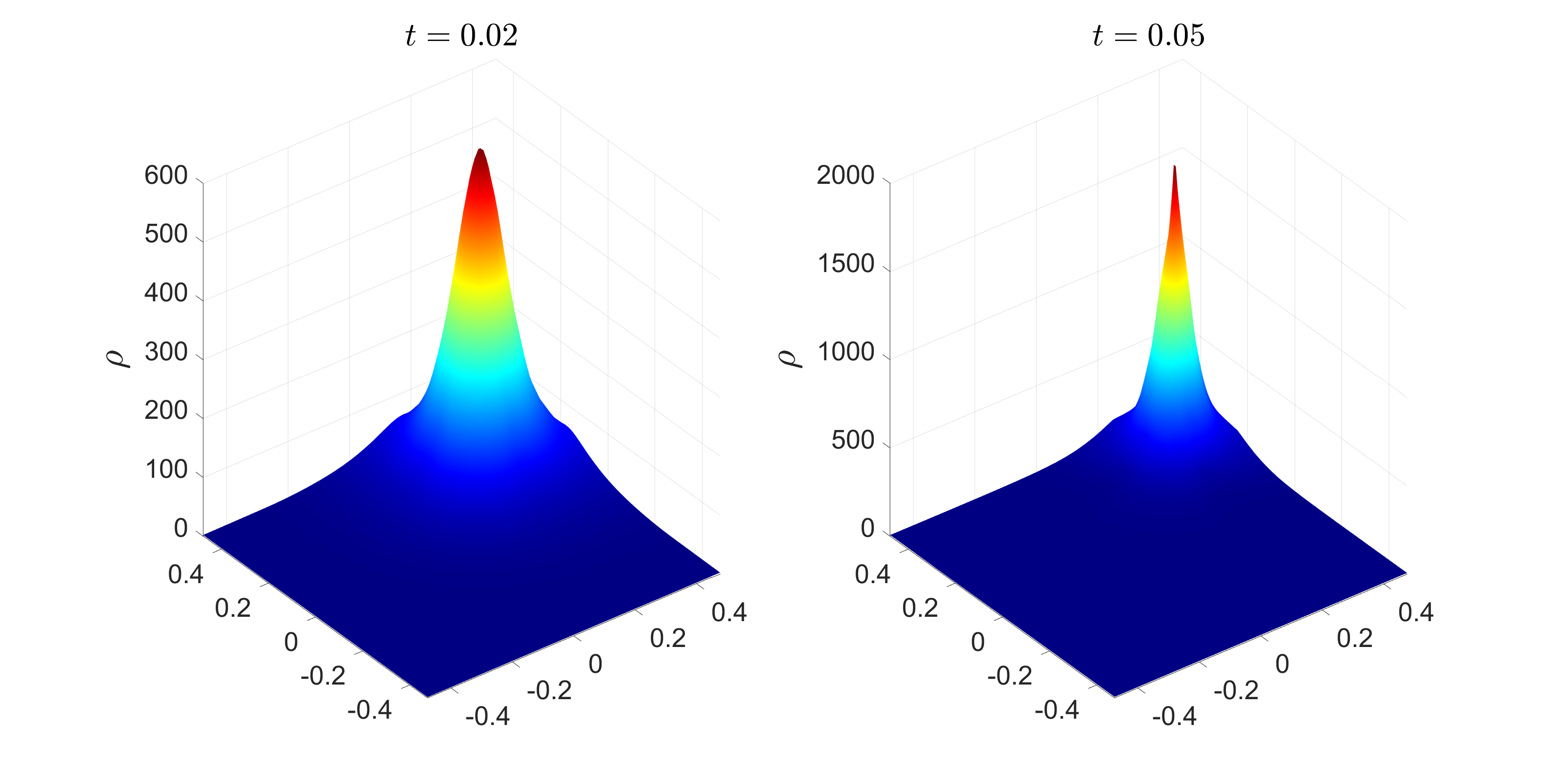}}
\vskip20pt
\centerline{\includegraphics[trim=0cm 1.9cm 0cm 0.7cm,clip,width=0.9\textwidth]{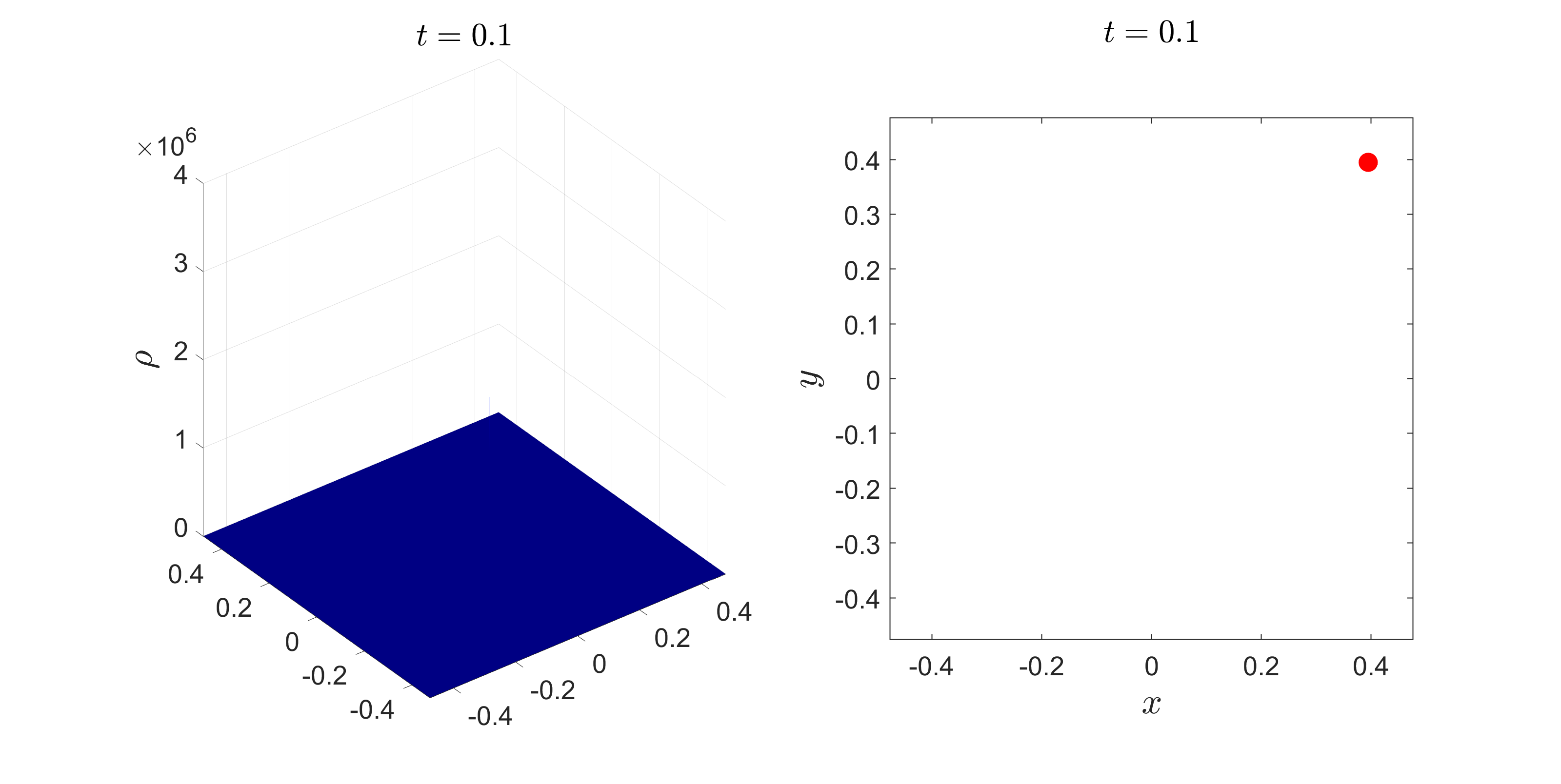}}
\caption{\sf Example 4: $\rho(\bm x,0.02)$ (upper left), $\rho(\bm x,0.05)$ (upper right), $\rho(\bm x,0.1)$ (lower left), and the location
of the ``dominating'' particle at $t=0.1$ (lower right), obtained using the hybrid FDP method with $\Delta=1/20$.\label{fig44}}
\end{figure}

We compute the solution by the proposed hybrid FDP method on a mesh of size $\Delta=1/20$ until the final time $t=0.1$. In Figure
\ref{fig44}, we plot the density $\rho$ at different times together with the location of the ``dominating'' particle at $t=0.1$ (by that
time the solution has blown up and most of the mass is concentrated at the ``dominating'' particle), which is, as one can clearly see,
inside the domain. This result is consistent with both the theoretical result in \cite{Herrero1997} and the numerical results reported in
\cite{Chertock2019a,AcostaSoba2023}.

\subsection{Parabolic-Elliptic Case ($\tau=0$)}\label{sec42}
In this section, we consider an extremely challenging numerical example introduced in \cite{Kurganov2014}.

\subsubsection*{Example 5---Two-Species Blowup at the Center of the Domain}
In this example, we consider the system \eref{KS} with $\tau=0$, $\gamma_1=\gamma_2=\zeta=\nu_1=\nu_2=\nu=\chi_1=1$, $\chi_2=20$, and
subject to the following initial conditions:
\begin{equation*}
\rho_1(x,y,0)=\rho_2(x,y,0)=50\,{\rm e}^{-100(x^2+y^2)},\quad(x,y)\in\Omega=[-1,1]\times[-1,1].
\end{equation*}
According to \cite{Espejo2012,Espejo2013}, both $\rho_1$ and $\rho_2$ are expected to blow up {\em simultaneously} in finite time. However,
as demonstrated in \cite{Kurganov2014,Chertock2018,Chertock2019}, $\rho_1$ and $\rho_2$ are expected to undergo different blowup patterns:
While $\rho_2$ develops a $\delta$-type singularity, $\rho_1$ grows up algebraically.

\begin{figure}[ht!]
\centerline{\includegraphics[trim=0cm 2.2cm 0cm 0.9cm,clip,width=0.9\textwidth]{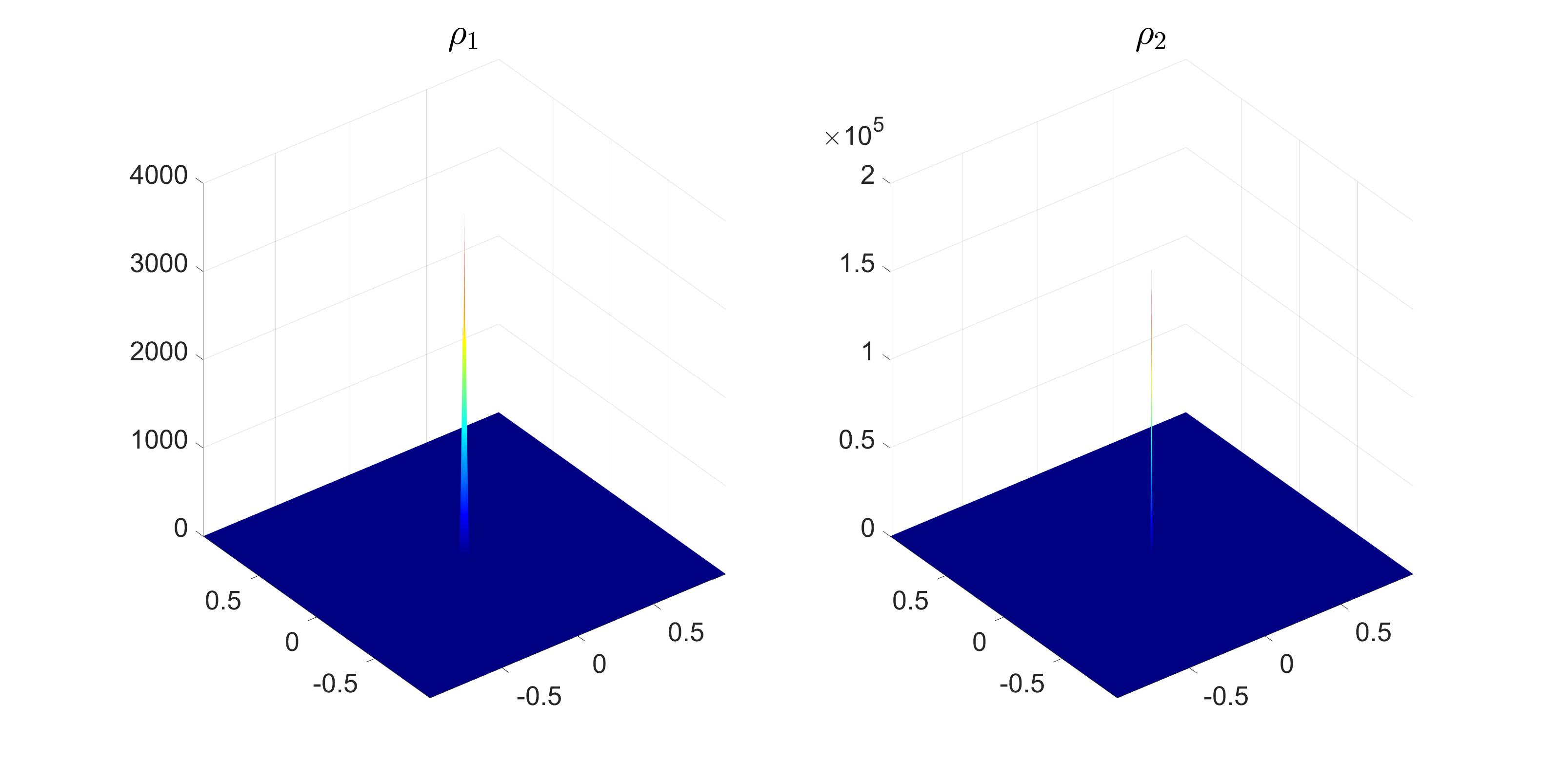}}
\vskip20pt
\centerline{\includegraphics[trim=0cm 2.2cm 0cm 0.9cm,clip,width=0.9\textwidth]{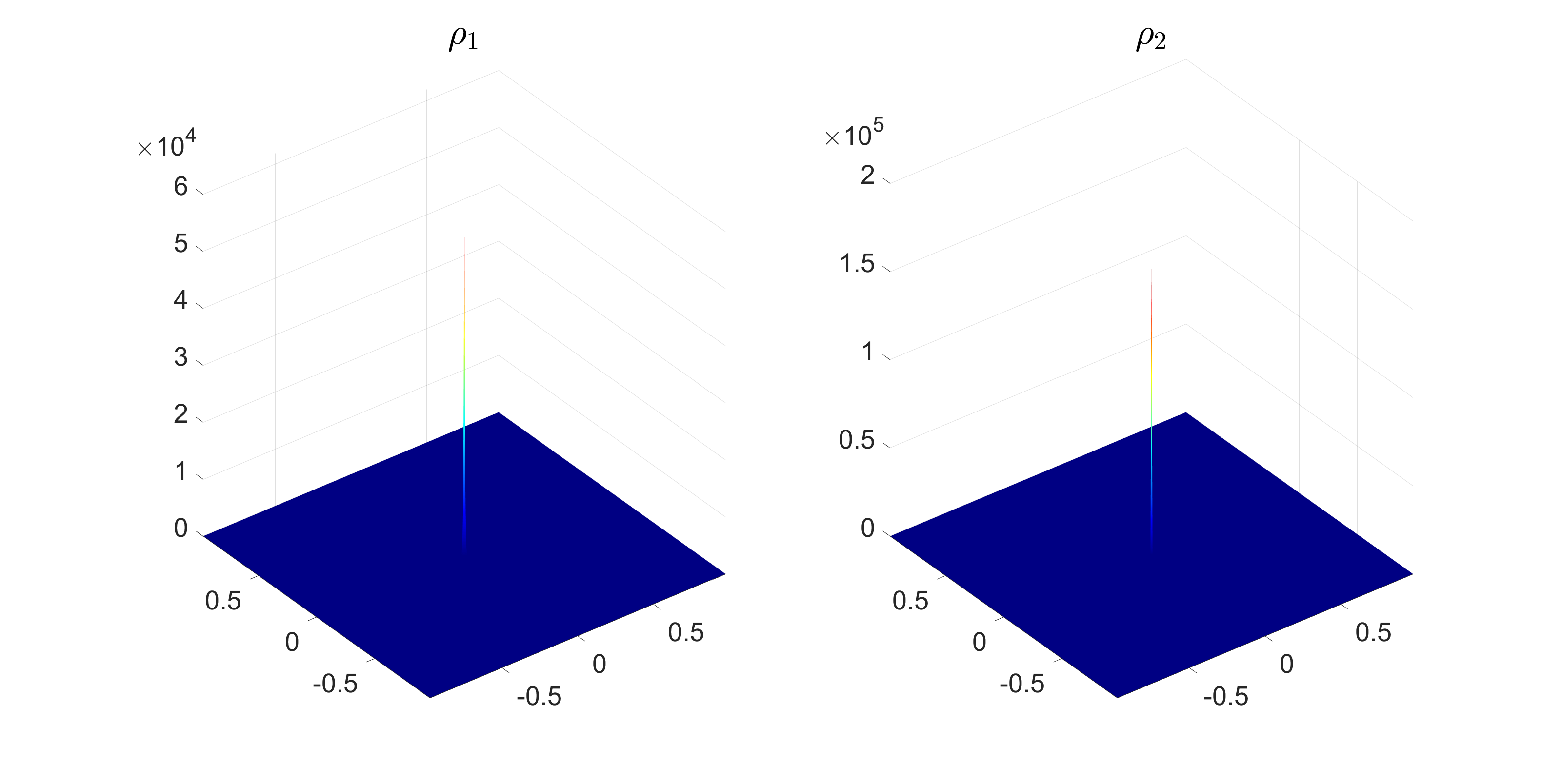}}
\caption{\sf Example 5: $\rho_1(\bm x,0.003)$ (upper left), $\rho_2(\bm x,0.003)$ (upper right), $\rho_1(\bm x,0.0033)$ (lower left), and
$\rho_2(\bm x,0.0033)$ (lower right), obtained using the hybrid FDP method with $\Delta=1/20$.\label{fig45}}
\end{figure}

We compute the solution until the same (as in \cite{Kurganov2014,Chertock2018,Chertock2019}) final time $t=0.0033$ using the proposed hybrid
FDP method on a mesh of size $1/20$. The densities $\rho_1$ and $\rho_2$, obtained at times $t=0.003$ and 0.0033 are shown in Figure
\ref{fig45}. We would like to emphasize that achieving a high resolution of blowup behavior in $\rho_1$ is a particularly challenging task.
Neither the second-order hybrid finite-volume-finite-difference \cite{Kurganov2014} nor its fourth-order version \cite{Chertock2018} were
able to clearly capture the algebraic blowup of $\rho_1$. Some improvement was achieved in \cite{Chertock2019}, where an AMM finite-volume
upwind method was developed and applied to the studied IBVP, but the maximum of the blowed up $\rho_1$ on the finest mesh there was about 75
(compare with the corresponding maximum of $\rho_2$, which was about $1.3\times10^5$). One can observe in Figure \ref{fig45}, our hybrid FDP
method leads to a considerably higher resolution of $\rho_1$ compared to the schemes in \cite{Kurganov2014,Chertock2018,Chertock2019}: The
maximum values of $\rho_1$ and $\rho_2$ at $t=0.0033$ are about $6.2\times10^4$ and $1.6\times 10^5$, respectively. The qualitative jump in
the resolution of $\rho_1$ is attributed to the low-dissipation nature of the sticky particle method.

\begin{figure}[ht!]
\centerline{\includegraphics[width=0.45\textwidth]{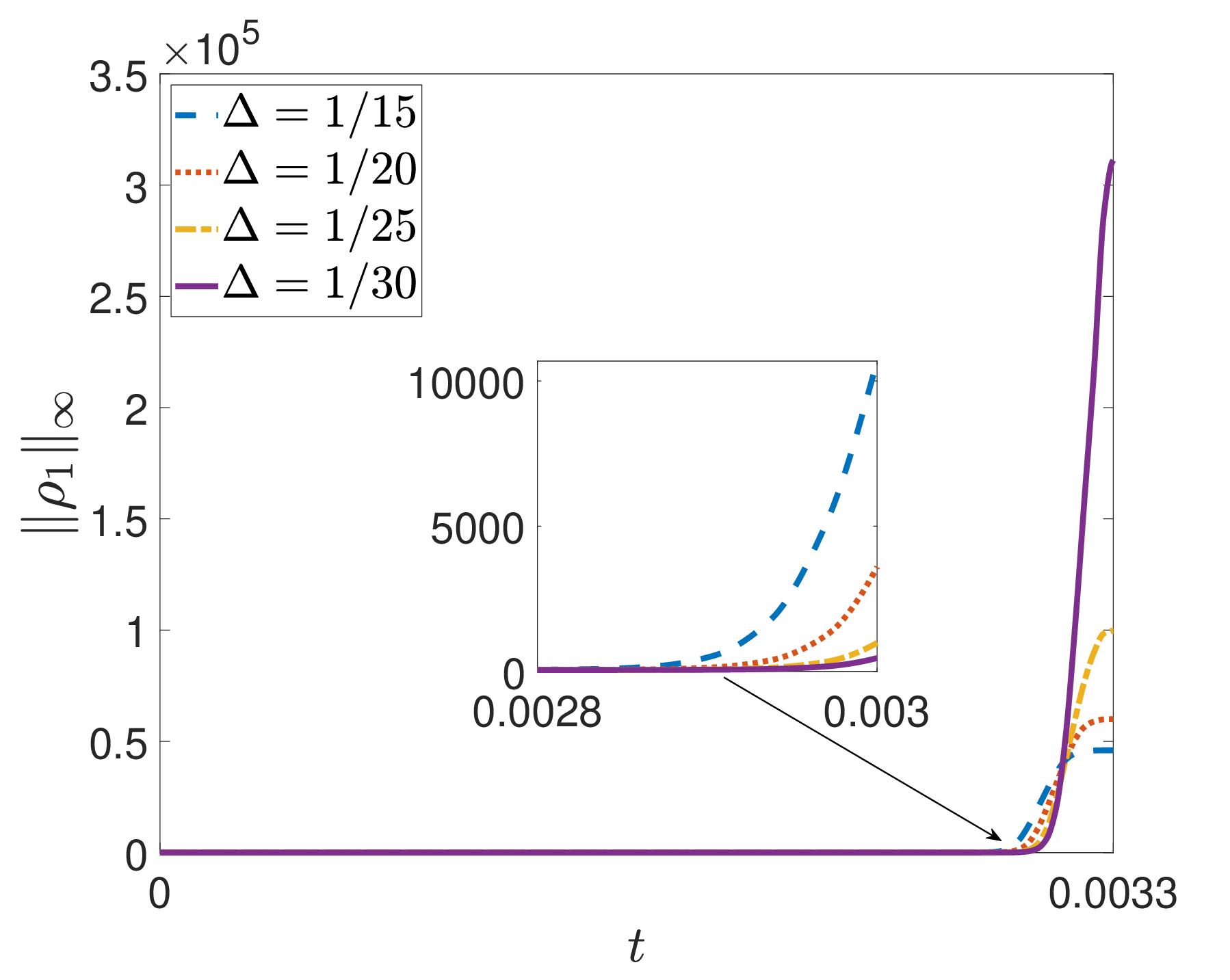}
            \includegraphics[width=0.45\textwidth]{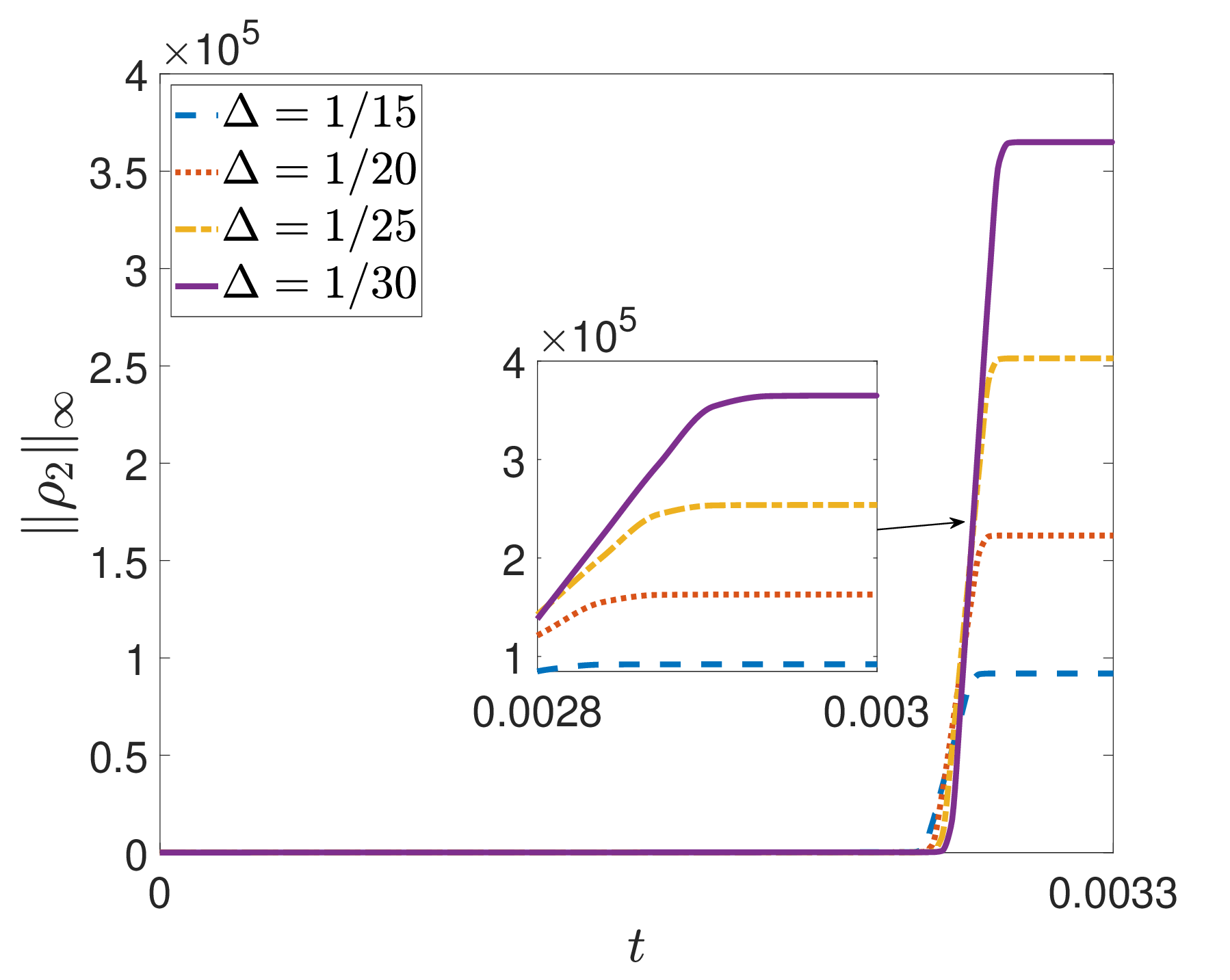}}
\caption{\sf Example 5: The maximums of $\rho_1$ (left) and $\rho_2$ (right) computed using $\Delta=1/15$, 1/20, 1/25, and 1/30.
\label{fig46}}
\end{figure}

In Figure \ref{fig46}, we present time evolution of the maximum values of $\rho_1$ and $\rho_2$ for $\Delta=1/15$, 1/20, 1/25, and 1/30. As
one can see, by refining the mesh from $\Delta=1/15$ to 1/30, the maximum of $\rho_1$ increases by a factor of about 6, which exhibits a
much faster increase compared those obtained in \cite{Kurganov2014,Chertock2018,Chertock2019}. This clearly demonstrates that the hybrid FDP
method outperforms its grid-based counterparts.

It should also be observed that the blowup time for the PKS-type system \eref{KS} cannot be obtained theoretically and therefore it is
important to determine it numerically. Figure \ref{fig46} suggests that the solution computed by the hybrid FDP method blows up by time
$t=0.00294$, at which the entire mass of the second species concentrates at one point (noticably, this time is the same for different
$\Delta$). We showcase this by zooming the maximum curves in Figure \ref{fig46} at the time interval $[0.028,0.003]$, which helps one to
better observe a difference between the blowup behavior of the first and second species. We note that the observed blowup time ($t=0.00294$)
is smaller the one ($t=0.0033$) reported in \cite{Chertock2019}; see also \cite{Kurganov2014,Chertock2018}. This is attributed to a
low-dissipation nature of the sticky particle method as the grid-based methods in \cite{Kurganov2014,Chertock2018,Chertock2019} contain
numerical dissipation, which causes a delay in blowup for both $\rho_1$ and $\rho_2$.

\section{Conclusion}\label{sec5}
In this paper, we have designed a novel hybrid finite-difference-particle (FDP) methods for the Patlak-Keller-Segel-type chemotaxis systems
in either parabolic-parabolic or parabolic-elliptic form. Our approach uses a sticky particle method for the cell density equation(s) and a
second-order finite-difference scheme for the chemoattractant equation. The performance of the proposed hybrid FDP methods have been
examined on a comprehensive series of challenging examples. Thanks to the low-dissipation nature of the sticky particle method, the proposed
scheme is particularly adept at capturing the blowing up solutions in both one- and two-species cases. The numerical results obtained using
the hybrid FDP methods exhibit superior resolution compared to traditional grid-based methods.

\section*{Acknowledgments}
The work of A. Chertock was supported in part by NSF grant DMS-2208438. The work of S. Cui was supported in part by Shenzhen Science and
Technology Program (grant No. RCJC20221008092757098). The work of A. Kurganov was supported in part by NSFC grant 12171226 and by the fund
of the Guangdong Provincial Key Laboratory of Computational Science and Material Design (No. 2019B030301001).


\bibliography{refs}

\begin{thebibliography}{10}

\bibitem{AcostaSoba2023}
{\sc D.~Acosta-Soba, F.~Guill\'{e}n-Gonz\'{a}lez, and J.~R.
  Rodr\'{\i}guez-Galv\'{a}n}, {\em An unconditionally energy stable and
  positive upwind {DG} scheme for the {K}eller-{S}egel model}, J. Sci. Comput.,
  97 (2023), pp.~Paper No. 18, 27.

\bibitem{BCK05}
{\sc M.~Bergdorf, G.-H. Cottet, and P.~Koumoutsakos}, {\em Multilevel adaptive
  particle methods for convection-diffusion equations}, Multiscale Model.
  Simul., 4 (2005), pp.~328--357.

\bibitem{BJ}
{\sc F.~Bouchut and F.~James}, {\em Duality solutions for pressureless gases,
  monotone scalar conservation laws, and uniqueness}, Comm. Partial
  Differential Equations, 24 (1999), pp.~2173--2189.

\bibitem{Chertock2017a}
{\sc A.~Chertock}, {\em A practical guide to deterministic particle methods},
  in Handbook of numerical methods for hyperbolic problems, vol.~18 of Handb.
  Numer. Anal., Elsevier/North-Holland, Amsterdam, 2017, pp.~177--202.

\bibitem{Chertock2017}
{\sc A.~Chertock, S.~Cui, and A.~Kurganov}, {\em Hybrid finite-volume-particle
  method for dusty gas flows}, SMAI J. Comput. Math., 3 (2017), pp.~139--180.

\bibitem{CDKK}
{\sc A.~Chertock, C.~R. Doering, E.~Kashdan, and A.~Kurganov}, {\em A fast
  explicit operator splitting method for passive scalar advection}, J. Sci.
  Comput., 45 (2010), pp.~200--214.

\bibitem{CDTM}
{\sc A.~Chertock, P.~Du~Toit, and J.~E. Marsden}, {\em Integration of the
  {EPD}iff equation by particle methods}, ESAIM Math. Model. Numer. Anal., 46
  (2012), pp.~515--534.

\bibitem{Chertock2018}
{\sc A.~Chertock, Y.~Epshteyn, H.~Hu, and A.~Kurganov}, {\em High-order
  positivity-preserving hybrid finite-volume-finite-difference methods for
  chemotaxis systems}, Adv. Comput. Math., 44 (2018), pp.~327--350.

\bibitem{CKpart}
{\sc A.~Chertock and A.~Kurganov}, {\em On a practical implementation of
  particle methods}, Appl. Numer. Math., 56 (2006), pp.~1418--1431.

\bibitem{Chertock2008}
\leavevmode\vrule height 2pt depth -1.6pt width 23pt, {\em A second-order
  positivity preserving central-upwind scheme for chemotaxis and haptotaxis
  models}, Numer. Math., 111 (2008), pp.~169--205.

\bibitem{CKL_2CH}
{\sc A.~Chertock, A.~Kurganov, and Y.~Liu}, {\em Finite-volume-particle methods
  for the two-component {C}amassa-{H}olm system}, Commun. Comput. Phys., 27
  (2020), pp.~480--502.

\bibitem{Chertock2019a}
{\sc A.~Chertock, A.~Kurganov, M.~Luk\'{a}\v{c}ov\'{a}-Medvi\v{d}ov\'{a}, and
  c.~N. \"{O}zcan}, {\em An asymptotic preserving scheme for kinetic chemotaxis
  models in two space dimensions}, Kinet. Relat. Models, 12 (2019),
  pp.~195--216.

\bibitem{Chertock2019}
{\sc A.~Chertock, A.~Kurganov, M.~Ricchiuto, and T.~Wu}, {\em Adaptive moving
  mesh upwind scheme for the two-species chemotaxis model}, Comput. Math.
  Appl., 77 (2019), pp.~3172--3185.

\bibitem{Chertock2007}
{\sc A.~Chertock, A.~Kurganov, and Y.~Rykov}, {\em A new sticky particle method
  for pressureless gas dynamics}, SIAM J. Numer. Anal., 45 (2007),
  pp.~2408--2441.

\bibitem{CLdisp}
{\sc A.~Chertock and D.~Levy}, {\em Particle methods for dispersive equations},
  J. Comput. Phys., 171 (2001), pp.~708--730.

\bibitem{CLP13}
{\sc A.~Chertock, J.-G. Liu, and T.~Pendleton}, {\em Elastic collisions among
  peakon solutions for the {C}amassa-{H}olm equation}, Appl. Numer. Math., 93
  (2015), pp.~30--46.

\bibitem{Cho}
{\sc A.~J. Chorin}, {\em Numerical study of slightly viscous flow}, J. Fluid
  Mech., 57 (1973), pp.~785--796.

\bibitem{Conca2011}
{\sc C.~Conca, E.~Espejo, and K.~Vilches}, {\em Remarks on the blowup and
  global existence for a two species chemotactic {K}eller-{S}egel system in
  {$\Bbb R^2$}}, European J. Appl. Math., 22 (2011), pp.~553--580.

\bibitem{CotKou}
{\sc G.-H. Cottet and P.~D. Koumoutsakos}, {\em Vortex methods. Theory and
  practice}, Cambridge University Press, Cambridge, 2000.

\bibitem{CKM13}
{\sc S.~Cui, A.~Kurganov, and A.~Medovikov}, {\em Particle methods for {PDE}s
  arising in financial modeling}, Appl. Numer. Math., 93 (2015), pp.~123--139.

\bibitem{Degond1989}
{\sc P.~Degond and S.~Mas-Gallic}, {\em The weighted particle method for
  convection-diffusion equations. {I}. {T}he case of an isotropic viscosity},
  Math. Comp., 53 (1989), pp.~485--507.

\bibitem{DMus}
{\sc P.~Degond and F.-J. Mustieles}, {\em A deterministic approximation of
  diffusion equations using particles}, SIAM J. Sci. Statist. Comput., 11
  (1990), pp.~293--310.

\bibitem{Eldredge2002}
{\sc J.~D. Eldredge, A.~Leonard, and T.~Colonius}, {\em A general deterministic
  treatment of derivatives in particle methods}, Journal of Computational
  Physics, 180 (2002), pp.~686--709.

\bibitem{Espejo2013}
{\sc E.~Espejo, K.~Vilches, and C.~Conca}, {\em Sharp condition for blow-up and
  global existence in a two species chemotactic {K}eller-{S}egel system in
  {$\Bbb{R}^2$}}, European J. Appl. Math., 24 (2013), pp.~297--313.

\bibitem{Espejo2012}
{\sc E.~E. Espejo, A.~Stevens, and T.~Suzuki}, {\em Simultaneous blowup and
  mass separation during collapse in an interacting system of chemotactic
  species}, Differential Integral Equations, 25 (2012), pp.~251--288.

\bibitem{GKS}
{\sc S.~Gottlieb, D.~Ketcheson, and C.-W. Shu}, {\em Strong stability
  preserving {R}unge-{K}utta and multistep time discretizations}, World
  Scientific Publishing Co. Pte. Ltd., Hackensack, NJ, 2011.

\bibitem{Gottlieb2001}
{\sc S.~Gottlieb, C.-W. Shu, and E.~Tadmor}, {\em Strong stability-preserving
  high-order time discretization methods}, SIAM Rev., 43 (2001), pp.~89--112.

\bibitem{Griebel00}
{\sc M.~Griebel and M.~A. Schweitzer}, {\em A particle-partition of unity
  method for the solution of elliptic, parabolic, and hyperbolic {PDE}s}, SIAM
  J. Sci. Comput., 22 (2000), pp.~853--890 (electronic).

\bibitem{Herrero1996}
{\sc M.~A. Herrero and J.~J.~L. Vel\'{a}zquez}, {\em Chemotactic collapse for
  the {K}eller-{S}egel model}, J. Math. Biol., 35 (1996), pp.~177--194.

\bibitem{Herrero1997}
\leavevmode\vrule height 2pt depth -1.6pt width 23pt, {\em A blow-up mechanism
  for a chemotaxis model}, Ann. Scuola Norm. Sup. Pisa Cl. Sci. (4), 24 (1997),
  pp.~633--683.

\bibitem{Hillen2009}
{\sc T.~Hillen and K.~J. Painter}, {\em A user's guide to {PDE} models for
  chemotaxis}, J. Math. Biol., 58 (2009), pp.~183--217.

\bibitem{Keller1970}
{\sc E.~F. Keller and L.~A. Segel}, {\em Initiation of slime mold aggregation
  viewed as an instability}, J. Theoret. Biol., 26 (1970), pp.~399--415.

\bibitem{KS711}
\leavevmode\vrule height 2pt depth -1.6pt width 23pt, {\em Model for
  chemotaxis}, J. Theor. Biol., 30 (1971), pp.~225--234.

\bibitem{Kurganov2014}
{\sc A.~Kurganov and M.~Luk\'{a}\v{c}ov\'{a}-Medvi\v{d}ov\'{a}}, {\em Numerical
  study of two-species chemotaxis models}, Discrete Contin. Dyn. Syst. Ser. B,
  19 (2014), pp.~131--152.

\bibitem{MGP}
{\sc S.~Mas-Gallic and F.~Poupaud}, {\em Approximation of the transport
  equation by a weighted particle method}, Transport Theory and Stat. Phys., 17
  (1988), pp.~311--345.

\bibitem{Othmer1997}
{\sc H.~G. Othmer and A.~Stevens}, {\em Aggregation, blowup, and collapse: the
  {ABC}s of taxis in reinforced random walks}, SIAM J. Appl. Math., 57 (1997),
  pp.~1044--1081.

\bibitem{Pat}
{\sc C.~S. Patlak}, {\em Random walk with persistence and external bias}, Bull.
  Math. Biophys., 15 (1953), pp.~311--338.

\bibitem{Pes}
{\sc C.~S. Peskin}, {\em The immersed boundary method}, Acta Numer., 11 (2002),
  pp.~479--517.

\bibitem{Puc}
{\sc E.~G. Puckett}, {\em Vortex methods: an introduction and survey of
  selected research topics}, in Incompressible computational fluid dynamics:
  trends and advances, Cambridge Univ. Press, Cambridge, 2008, pp.~335--407.

\bibitem{Rav}
{\sc P.-A. Raviart}, {\em An analysis of particle methods}, in Numerical
  methods in fluid dynamics (Como, 1983), vol.~1127 of Lecture Notes in Math.,
  Springer, Berlin, 1985, pp.~243--324.

\bibitem{Senba2001}
{\sc T.~Senba and T.~Suzuki}, {\em Parabolic system of chemotaxis: blowup in a
  finite and the infinite time}, Methods Appl. Anal., 8 (2001), pp.~349--367.
\newblock IMS Workshop on Reaction-Diffusion Systems (Shatin, 1999).

\bibitem{Wolansky2002}
{\sc G.~Wolansky}, {\em Multi-components chemotactic system in the absence of
  conflicts}, European J. Appl. Math., 13 (2002), pp.~641--661.

\end{thebibliography}
\bibliographystyle{siam}

\end{document}